\documentclass[a4paper]{amsart}
\usepackage{amssymb}
\usepackage{graphicx}
\usepackage{psfrag}
\usepackage[all]{xy}
\xyoption{matrix}
\xyoption{arrow}

\begin{document}

\newtheorem{theo}{Theorem}[section]
\newcommand{\bthm}{\begin{theo}}
\newcommand{\ethm}{\end{theo}}
\newtheorem{coro}[theo]{Corollary}
\newcommand{\bcor}{\begin{coro}}
\newcommand{\ecor}{\end{coro}}
\newtheorem{prop}[theo]{Proposition}
\newcommand{\bprop}{\begin{prop}}
\newcommand{\eprop}{\end{prop}}
\newtheorem{claim}[theo]{Lemma}
\newcommand{\blemma}{\begin{claim}}
\newcommand{\elemma}{\end{claim}}
\newtheorem{conjec}[theo]{Conjecture}
\newcommand{\bconj}{\begin{conjec}}
\newcommand{\econj}{\end{conjec}}

\theoremstyle{definition}
\newtheorem{defo}[theo]{Definition}
\newcommand{\bdfn}{\begin{defo}}
\newcommand{\edfn}{\end{defo}}
\newtheorem{rmk}[theo]{Remark}
\newcommand{\brk}{\begin{rmk}}
\newcommand{\erk}{\end{rmk}}
\newtheorem{exple}[theo]{Example}
\newcommand{\bex}{\begin{exple}}
\newcommand{\eex}{\end{exple}}

\newcommand{\note}{\mbox{\textbf{Note: \hspace{0.2cm}}}}


\title[Explicit construction of companion bases]{Explicit construction of companion bases}

\author[Parsons]{Mark James Parsons}
\address{
Institut f\"{u}r Mathematik und Wissenschaftliches Rechnen \\
Universit\"{a}t Graz \\
Heinrichstra\ss e 36 \\
A-8010 Graz \\
Austria
}
\email{mark.parsons@uni-graz.at}


\begin{abstract}
We present a simple procedure for explicitly constructing a companion basis for a quiver mutation equivalent to a simply-laced Dynkin quiver.
\end{abstract}

\date{2 December 2013}



\thanks{This work was supported by the
\textit{Engineering and Physical Sciences Research Council} and by the \textit{Austrian Science Fund (FWF): Project No. P25141-N26}.}

\maketitle


\section*{Introduction}

The cluster algebras of finite type were shown to have a classification by Dynkin diagrams in \cite{FZ}. We consider the cluster algebra $\mathcal{A}$ associated to a simply-laced Dynkin diagram $\Delta$. The exchange matrices of $\mathcal{A}$ are skew-symmetric integer matrices with entries not exceeding 1 in absolute value and can therefore be represented as quivers, which we shall refer to as the quivers of mutation type $\Delta$. Fix such a quiver $\Gamma$ and denote its set of vertices by $\Gamma_0$. In \cite{Par2} (originally in \cite{Par1}), motivated by work in \cite{BGZ}, we introduced the notion of a companion basis for $\Gamma$ as a certain subset of the corresponding root system $\Phi$ of type $\Delta$. Specifically, a \textit{companion basis} for $\Gamma$ is a subset $\{ \gamma_x \colon x \in \Gamma_0 \} \subseteq \Phi$ whose elements form a $\mathbb{Z}$-basis for the integral root lattice $\mathbb{Z} \Phi$ and such that for distinct vertices $x,y \in \Gamma_0$, $\left | (\gamma_x , \gamma_y) \right |$ is equal to the number of edges joining $x$ and $y$ in the underlying unoriented graph for $\Gamma$.

This article continues the study of companion bases initiated in \cite{Par2}. The focus here is on producing explicit companion bases for quivers. In \cite{Par2}, we already presented a method for finding a companion basis for any given quiver of simply-laced Dynkin mutation type. Indeed, we noted there that a simple system of a root system of simply-laced Dynkin type is a companion basis for a Dynkin quiver (that is, a quiver whose underlying graph is a simply-laced Dynkin diagram) of the same type. Furthermore, we introduced a companion basis mutation procedure that, given a companion basis for a quiver of simply-laced Dynkin mutation type, produces a companion basis for any mutation of that quiver. While this, at least in theory, enables us to find a companion basis for any quiver of simply-laced Dynkin mutation type, it has a major drawback that can make it difficult to apply in practice. Namely, that it requires us to find a sequence of quiver mutations taking us from a quiver for which we already have a companion basis to the quiver for which we desire to produce one. Here, we show that this problem can be avoided by outlining a method for explicitly constructing a companion basis for a given quiver of simply-laced Dynkin mutation type, which depends only on the structure and properties of the quiver. We do this initially for the type $A$ case, before showing how the approach used in this case can be used and adapted to deal also with the type $D$ case. With the aim of keeping the paper to a reasonable length, albeit at the expense of completeness, we refrain from writing down explicit companion bases for each of the quivers arising in types $E_6$, $E_7$ and $E_8$. Since the problem is finite in these types, it is fairly easy, although time-consuming, to compute companion bases case-by-case using the method explained above.

Our companion basis construction procedure consists of two steps. The first is to apply a labelling to the vertices of the quiver for which we wish to construct a companion basis. This then enables us to read-off (that is, to directly write down) a companion basis for that quiver. Fix a simple system $\Pi = \{ \alpha_1 , \ldots , \alpha_n \}$ of $\Phi$. The procedure for labelling a quiver of mutation type $A$ is presented in Section~\ref{labelling section}. Given a labelled quiver of mutation type $A$, Theorem~\ref{companion basis type A} (which uses terminology from Definition~\ref{3-cycle vertex types}) then provides us with a companion basis for that quiver. The quivers of mutation type $D$ can be separated into four types, as described by Vatne in \cite{Vat}. We recall this description at the beginning of Section~\ref{type D construction}. We then consider each type in turn, first explaining the labelling procedure for a given quiver of that type, before establishing a result which provides us with a companion basis for the labelled quiver. The results giving companion bases for labelled quivers of types I, II and III are Theorems~\ref{type D type I}, \ref{type D type II} and \ref{type D type III} respectively. For convenience, we split type IV into two cases. Theorems~\ref{type D type IV case 1} and \ref{type D type IV case 2} provide companion bases in these cases. (Refer to Definition~\ref{vertex types in type D} for the terminology used in Theorems~\ref{type D type I} to \ref{type D type IV case 2}.)


\section{Motivation and applications}
\label{applications section}

The initial motivation for the definition and study of companion bases arose from work in \cite{BGZ} on the recognition of cluster algebras of finite type, in which positive quasi-Cartan matrices were considered. We recall that a \textit{quasi-Cartan matrix} is a symmetrizable matrix for which each diagonal entry is 2, and that such a matrix is said to be \textit{positive} if its principal minors are all positive. Their main result, \cite[Theorem 1.2]{BGZ}, gave a necessary and sufficient condition for a given skew-symmetrizable integer matrix to give rise to a cluster algebra of finite type. Part of this condition is for the matrix to have a positive quasi-Cartan companion. That is, for there to exist a positive quasi-Cartan matrix whose off-diagonal entries agree, up to sign, with those of the given matrix.

Let $B$ be the exchange matrix associated to $\Gamma$. We note that a companion basis for $\Gamma$ was originally defined \cite[Definition 4.1]{Par2} to be a $\mathbb{Z}$-basis $\{ \gamma_x \colon x \in \Gamma_0 \} \subseteq \Phi$ for $\mathbb{Z} \Phi$ such that the matrix $A = (a_{xy})$ given by $a_{xy} = (\gamma_x , \gamma_y)$ for all $x,y \in \Gamma_0$ is a positive quasi-Cartan companion of $B$. That this is equivalent to the definition we have given here follows immediately upon noting that the matrix of inner products of any $\mathbb{Z}$-basis of roots for $\mathbb{Z} \Phi$ is necessarily a positive quasi-Cartan matrix.

In view of \cite[Theorem 4.10 \& Corollary 4.11]{Par2}, the results we present in this article complete a description of the companion bases for a given quiver of simply-laced Dynkin mutation type. Indeed, \cite[Theorem 4.10 \& Corollary 4.11]{Par2} gives a complete description of the relationship between different companion bases for such a quiver, showing how all of the companion bases for a quiver can be expressed in terms of any given companion basis for that quiver. Consequently, from a single companion basis for a quiver, we can write down all of the companion bases for that quiver.

Another application of our companion basis construction procedure arises from \cite{BM} where it is shown how a given quiver of simply-laced Dynkin mutation type can be used to give a presentation of the associated finite crystallographic reflection group (Weyl group). We remark that by \cite[Theorem 6.8]{BM}, the reflections corresponding to the elements of any companion basis for the given quiver can be regarded as the generators in this presentation. Thus the results we obtain here give an explicit description of the generators appearing in the presentations of \cite{BM} in types $A$ and $D$. (We note that the results of \cite{BM} hold in finite type, also covering the non-simply-laced case.)

Other important factors motivating the study of companion bases in general are their links to cluster-tilted algebras, and moreover, the strong relationship between cluster-tilted algebras and cluster algebras. It was shown in \cite[Section 6]{BMR2} and \cite[Theorem 3.1]{CCS2} that the quivers of simply-laced Dynkin mutation type are precisely the quivers of the cluster-tilted algebras of simply-laced Dynkin type.

Let $\Lambda$ be the cluster-tilted algebra with quiver $\Gamma$. We showed in \cite[Theorem 5.3]{Par2} that, in Dynkin type $A$, the dimension vectors of the finitely generated indecomposable modules over $\Lambda$ can be obtained from any companion basis for $\Gamma$. In a result analogous to Gabriel's Theorem, they are obtained by expressing each positive root in terms of the companion basis, but replacing each integer with its absolute value in the corresponding coefficient tuple. Furthermore, we conjectured that this result may be extended to simply-laced Dynkin type \cite[Conjecture 6.3]{Par2}. That we could then use explicit companion bases to compute dimension vectors (as already established for type $A$) would provide us with a nice application of our companion basis construction procedure.

Independent work of Ringel \cite{Rin} gives an alternative approach for obtaining the dimension vectors of the finitely generated indecomposable modules over a cluster-tilted algebra, whenever the associated cluster-tilting object corresponds to a preprojective tilting module (the terminology he uses for such algebras is \textit{cluster-concealed}). This, in particular, includes the cluster-tilted algebras of simply-laced Dynkin type. We expect that the aforementioned conjecture could be established by understanding how companion bases fit into Ringel's set-up. We note that it is enough to establish the conjecture for a single companion basis for each quiver of simply-laced Dynkin mutation type due to \cite[Proposition 6.2]{Par2}.

Dimension vectors of finitely generated indecomposable modules over cluster-tilted algebras of simply-laced Dynkin type play an important role in the relationship between these algebras and the associated cluster algebras. Denote the set of dimension vectors of the finitely generated indecomposable modules over $\Lambda$ by $\mathrm{Dim}(\Lambda)$. Let $\mathcal{A}_{\bullet}(B)$ be the cluster algebra (of simply-laced Dynkin type $\Delta$) with principal coefficients whose initial exchange matrix is $B$. It was shown in \cite[Theorem 4.4 \& Remark 4.5]{CCS2}, \cite[Theorem 2.2]{BMR1} that the set of non-initial $d$-vectors of $\mathcal{A}_{\bullet}(B)$ coincides with $\mathrm{Dim}(\Lambda)$. In addition, it was established in \cite[Theorem 6]{NCh} (see also \cite[Section 9]{DWZ}) that the set of positive $c$-vectors of $\mathcal{A}_{\bullet}(B)$ also coincides with these sets. For further information on these interesting connections, we also refer the reader to \cite{NS} which gives explicit diagrammatic descriptions of the $c$- and $d$-vectors of a cluster algebra of finite type with any initial exchange matrix.

It is hoped that this article, in providing explicit companion bases for quivers of simply-laced Dynkin mutation type, will facilitate further research into the links between companion bases and both cluster algebras and cluster-tilted algebras.

We conclude this section with a brief comment on the non-simply-laced case.

The definition of companion bases can be naturally extended to the finite type case. In this setting, since the exchange matrices appearing in cluster algebras are skew-symmetrizable (rather than skew-symmetric), quivers are replaced by valued quivers. In the shortly forthcoming article \cite{Par3}, the results of the present article are extended to the finite type case. Explicit companion basis construction procedures for valued quivers of mutation types $B$ and $C$ are obtained by again making use of the procedure presented here for quivers of mutation type $A$. Companion bases for valued quivers of types $F_4$ and $G_2$ are computed directly, using companion basis mutation.


\section{Labelling procedure for quivers arising in Dynkin type $A_n$}
\label{labelling section}

In this section, we focus on the Dynkin type $A$ case. In particular, we introduce a procedure for labelling the vertices of the quiver $\Gamma$. This labelling will subsequently enable us to directly write down (read off) a companion basis for $\Gamma$. We first need to explain the structure of $\Gamma$. This can be done through consideration of the triangulations of a regular $(n+3)$-gon $\mathbb{P}_{n+3}$. To a triangulation $\mathbb{T}$ of $\mathbb{P}_{n+3}$, we can associate a connected quiver $Q_{\mathbb{T}}$ as in \cite[Section 2.4]{CCS1}: The vertices of $Q_{\mathbb{T}}$ correspond to the diagonals in $\mathbb{T}$. Let $i$ and $j$ be vertices in $Q_{\mathbb{T}}$ corresponding to diagonals $d_i$ and $d_j$ respectively. Then, there is an arrow from $i$ to $j$ in $Q_{\mathbb{T}}$ if $d_i$ and $d_j$ bound a common triangle, and the angle of minimal rotation about the common point of $d_i$ and $d_j$ taking the line through $d_i$ to the line through $d_j$ is in the anticlockwise direction. It follows from \cite[Section 12]{FZ} that the quivers arising from the triangulations of $\mathbb{P}_{n+3}$ are precisely the quivers of mutation type $A_n$. In particular, we have that $\Gamma$ may be identified with the quiver $Q_{\mathbb{T}}$ arising from some triangulation $\mathbb{T}$ of $\mathbb{P}_{n+3}$.

From the triangulation $\mathbb{T}$, the structure of $\Gamma$ may be easily understood. The triangles appearing in $\mathbb{T}$ are of the following three types:

(I) Triangles that consist of one diagonal and two boundary edges of $\mathbb{P}_{n+3}$.

(II) Triangles that consist of two diagonals and one boundary edge of $\mathbb{P}_{n+3}$.

(III) Triangles that consist of three diagonals of $\mathbb{P}_{n+3}$.

\noindent We have that $\Gamma$ is connected and made up of oriented 3-cycles arising from triangles of type (III), joined together by linear sections (that is, full subquivers whose underlying graphs are Dynkin diagrams of type $A$) in such a way that no new cycles are introduced in the underlying (unoriented) graph. That is, in the underlying graph of $\Gamma$, the only cycles are 3-cycles arising from triangles of type (III). Moreover, for $n \geq 2$, we have that vertices have valency four if and only if they lie at a point where two 3-cycles meet, valency three if and only if they lie at a point where a linear section meets a 3-cycle, valency one if and only if they lie at the end of a linear section (not meeting a 3-cycle), and valency two otherwise. In the $n=1$ case, $\Gamma$ consists of a solitary vertex of valency zero. (This description of the quivers arising in Dynkin type $A_n$ was given independently in \cite{Par1, Par2} and \cite{BV}.)

Armed with the knowledge of the structure of $\Gamma$, we are able to complete the description of the cluster-tilted algebra $\Lambda$ with quiver $\Gamma$ (and hence of all cluster-tilted algebras of type $A_n$) by specifying the relations. Let $I_{\Gamma}$ be the admissible ideal of the path algebra $k\Gamma$ generated by the set of all paths in $\Gamma$ consisting of two consecutive arrows in any given 3-cycle, then we have from \cite[Theorem 4.2]{BMR3}, \cite[Theorem 4.1]{CCS2} that the cluster-tilted algebra $\Lambda$ is isomorphic to $\frac{k\Gamma}{I_{\Gamma}}$.

We now introduce the terminology we will require for our labelling procedure. Firstly, it will be useful for us to distinguish certain vertices in $\Gamma$.

\bdfn
\label{triangle vertex}
We call any vertex of $\Gamma$ belonging to a cyclically oriented 3-cycle a \textit{3-cycle vertex}.
\edfn

\bdfn
\label{end vertices}
A vertex of $\Gamma$ is said to be an \textit{end vertex} if it has valency zero, valency one, or is a 3-cycle vertex of valency two.
\edfn

So, the end vertices of $\Gamma$ are precisely those vertices corresponding to diagonals in $\mathbb{T}$ which bound a triangle in $\mathbb{T}$ of type (I).

Now, it is well known that any triangulation of a regular polygon with at least four sides must contain at least two triangles of type (I). (In general, a triangulation containing $m$ triangles of type (III) must have $m+2$ triangles of type (I).) Also, no diagonal in a triangulation of a regular polygon with at least five sides can bound two triangles of type (I). In particular, any quiver associated to a triangulation of a regular polygon with at least five sides must have at least two end vertices. We note also that any triangulation of a regular 4-gon must consist of a single diagonal bounding two triangles of type (I).

Important for the labelling procedure is the notion of a string in $\Gamma$. A \textit{string} in $\Gamma$ is a reduced walk in the quiver $\Gamma$ which avoids the zero relations and hence does not contain two consecutive arrows of any given 3-cycle in $\Gamma$. We recall that the cluster-tilted algebras of Dynkin type $A_n$ are known to be string algebras and that the isomorphism classes of the finitely generated indecomposable modules over such an algebra may be described in terms of the strings in its quiver (refer to \cite[Section 3]{BR} for both the definition of a string algebra and a description of the finitely generated indecomposable modules over such an algebra).

Let $i,j \in \Gamma_0$. Due to the structure of $\Gamma$, it is easily seen that there is a unique string $p$ from $i$ to $j$, which we call a \textit{trivial string} in the case where $i = j$. By definition, $p$ does not pass through two consecutive arrows of any given 3-cycle in $\Gamma$, and none of the arrows appearing in $p$ appear more than once. Likewise, none of the vertices appearing in $p$ appear more than once.

Starting from $i$ and moving along $p$ towards $j$, we may pass through a number of 3-cycle vertices. We make the following definition.

\bdfn
\label{triangle vertex rel to p}
For any 3-cycle with two vertices appearing in $p$, we call the first vertex of that 3-cycle appearing in $p$ (when moving from $i$ towards $j$) a \textit{primary vertex relative to $p$}. We call the second vertex of that 3-cycle appearing in $p$ a \textit{secondary vertex relative to $p$}, and we call the vertex of that 3-cycle not appearing in $p$ a \textit{complementary vertex relative to $p$}.

Furthermore, if $x$, $y$ and $z$ are the three vertices of some 3-cycle in $\Gamma$ and are respectively primary, secondary and complementary vertices relative to $p$, then we call $y$ the \textit{secondary vertex corresponding to $x$ relative to $p$}, and we call $z$ the \textit{complementary vertex corresponding to $x$ relative to $p$}.
\edfn

It is worth noting that since two 3-cycles can meet in a vertex in $\Gamma$, then it is possible for a vertex to be a secondary vertex relative to $p$ with respect to one 3-cycle, and a primary vertex relative to $p$ with respect to another 3-cycle.

\bdfn
\label{section of quiver}
Suppose that $x$, $y$ and $z$ are the three vertices of a 3-cycle in $\Gamma$. We define the \textit{subquiver of $\Gamma$ rooted at $x$} to be the full subquiver of $\Gamma$ on all vertices that can be reached on strings starting at $x$ which do not pass through $y$ or $z$.
\edfn

\brk
\label{section remark}
For any vertex $v$ belonging to a 3-cycle in $\Gamma$, it is clear that the subquiver of $\Gamma$ rooted at $v$ arises as the quiver associated to a triangulation of a regular $(m+3)$-gon for some $m<n$ (where $m$ is the number of vertices in that subquiver). In addition, with $x$, $y$ and $z$ as in Definition~\ref{section of quiver} above, it is a simple observation that the subquivers of $\Gamma$ rooted at $x$, $y$ and $z$ are pairwise disjoint.
\erk

With all of the groundwork finally in place, we will now outline our procedure for labelling the vertices of a quiver of mutation type $A$. Each vertex will be labelled with a distinct natural number, with the lowest label being 1 and the highest being the number of vertices of the quiver. In order to label the vertices of a quiver, we must first choose an ordered pair of end vertices in that quiver, distinct if possible. Note that we have already observed that any quiver of mutation type $A$ with at least two vertices must have at least two distinct end vertices. The procedure for labelling a quiver with a given initial choice of an ordered pair of end vertices is an inductive one.

We first consider the quiver (of mutation type $A_1$) with a single vertex and no arrows. In this case, we label the vertex 1 whereby the labelling is completed.

For induction, we may suppose that we have obtained a labelling of the vertices of every quiver of mutation type $A$ with fewer than $n$ vertices, for any given initial choice of an ordered pair of end vertices (distinct if possible, as noted above).

Supposing now that $n \geq 2$, we consider the quiver $\Gamma$, and start by choosing an ordered pair of distinct end vertices in $\Gamma$. Label the first vertex of this ordered pair 1, and consider the (unique) string $p$ in $\Gamma$ from 1 to the other chosen end vertex. Starting from 1, move along $p$, labelling subsequent vertices consecutively $2,3,4, \ldots ,$ up to and including the first primary vertex $i$ relative to $p$. Note that 1 could be a primary vertex relative to $p$. Note also that there may be no primary vertices relative to $p$, in which case, the labelling procedure is completed here.

Denote the subquiver of $\Gamma$ rooted at the complementary vertex corresponding to $i$ relative to $p$ by $\Gamma'$. As a consequence of Remark~\ref{section remark}, it follows that $\Gamma'$ is a quiver of a mutation type $A$. Suppose that there are $a$ vertices in $\Gamma'$. We have that the complementary vertex corresponding to $i$ relative to $p$ is an end vertex in $\Gamma'$. We then obtain an ordered pair of end vertices in $\Gamma'$ by choosing another end vertex in $\Gamma'$, distinct if $a>1$, and setting this chosen end vertex to be the first vertex in the pair. By induction, we have a labelling of the vertices of $\Gamma'$ using the labels $1$ up to $a$, given this choice of an ordered pair of end vertices in $\Gamma'$. Add $i$ to each of the vertex labels in this labelling for $\Gamma'$, and then assign the labels thus obtained to the corresponding vertices in $\Gamma$.

Label the secondary vertex corresponding to $i$ relative to $p$ with $i+a+1$. Then, from $i+a+1$, continue along $p$ labelling subsequent vertices consecutively $i+a+2,i+a+3,\ldots ,$ and proceed as above for each subsequent primary vertex relative to $p$. The second of our initially chosen end vertices of $\Gamma$ will be the final vertex to be labelled, and will be labelled $n$.

\brk
\label{end vertex choices remark}
To label $\Gamma$ according to the above procedure, we initially choose a pair of end vertices in $\Gamma$. We also make a further choice of an end vertex for each primary vertex relative to some string considered in the labelling procedure. Because of these choices, there are potentially many different labellings of $\Gamma$ that can be obtained using the outlined procedure. Any such labelling can be used in the construction of a companion basis for $\Gamma$. Using a different labelling usually has the effect of changing the companion basis that is produced.
\erk

In order to illustrate the labelling procedure, we conclude this section by giving a detailed example of it in action.

\bex \rm
\label{labelling example}
Let $\Omega$ be the following quiver of mutation type $A_{11}$.
$$
\xymatrix @R=4pt @C=4pt {
{\Omega} & & & {\bullet} \ar[dd] \\
\\
 & & & {\bullet} \ar[dd] \\
\\
 & & & {\bullet} \ar[drr] \\
 & & & & & {\bullet} \ar[dll] & & {\bullet} \ar[ll] \\
 & & & {\bullet} \ar[uu] \ar[dd] \\
\\
 & & & {\bullet} \ar[ddr] \\
\\
{\bullet} \ar[rr]_<{a} & & {\bullet} \ar[uur] & & {\bullet} \ar[ll] \ar[rr]_>{b} & & {\bullet}
}
$$

We will use the prescribed labelling procedure to obtain a labelling of $\Omega$. The first step is to choose an ordered pair $(a,b)$ of end vertices in $\Omega$ (as shown above). We label the first vertex in this pair 1, and then consider the string $p$ in $\Omega$ from 1 to $b$. Starting from 1 and moving along $p$, labelling subsequent vertices in increments of one, we have that 2 is the first primary vertex relative to $p$.
$$
\xymatrix @R=4pt @C=4pt {
{\Omega} & & & {\bullet} \ar[dd] \\
\\
 & & & {\bullet} \ar[dd] \\
\\
 & & & {\bullet} \ar[drr] \\
 & & & & & {\bullet} \ar[dll] & & {\bullet} \ar[ll] \\
 & & & {\bullet} \ar[uu] \ar[dd]_>{c} \\
\\
 & & & {\bullet} \ar[ddr] \\
\\
{1} \ar[rr] & & {2} \ar[uur] & & {\bullet} \ar[ll] \ar[rr]_>{b} & & {\bullet} \\
\\
{} \ar[rrrrrr]_p & & & & & & {}
}
$$

We must now consider the subquiver $\Omega'$ of $\Omega$ rooted at the complementary vertex $c$ corresponding to 2 relative to $p$. Note that $\Omega'$ is a quiver of mutation type $A_7$.
$$
\xymatrix @R=4pt @C=4pt {
{\Omega'} & & & {\bullet} \ar[dd]_<{d} \\
\\
 & & & {\bullet} \ar[dd] \\
\\
 & & & {\bullet} \ar[drr] \\
 & & & & & {\bullet} \ar[dll] & & {\bullet} \ar[ll] \\
 & & & {\bullet} \ar[uu] \ar[dd]_>{c} \\
\\
 & & & {\bullet}
}
$$

\noindent Since $c$ is an end vertex in $\Omega'$, by choosing another end vertex $d$ in $\Omega'$, we obtain an ordered pair $(d,c)$ of end vertices. We now start to label the vertices of $\Omega'$, proceeding in the same manner as above. We label the first vertex of our ordered pair 1, and consider the string $p'$ in $\Omega'$ from 1 to $c$. Starting from 1 and moving along $p'$, labelling subsequent vertices in increments of one, we have that 3 is the first primary vertex relative to $p'$.
$$
\xymatrix @R=4pt @C=4pt {
{\Omega'} & {} \ar[dddddddd]_{p'} & & {1} \ar[dd] \\
\\
 & & & {2} \ar[dd] \\
\\
 & & & {3} \ar[drr] \\
 & & & & & {\bullet} \ar[dll] & & {\bullet} \ar[ll]^>{e} \\
 & & & {\bullet} \ar[uu] \ar[dd]_>{c} \\
\\
 & {} & & {\bullet}
}
$$

We must now consider the subquiver $\Omega''$ of $\Omega'$ rooted at the complementary vertex $e$ corresponding to 3 relative to $p'$.
$$
\xymatrix @R=4pt @C=4pt {
{\Omega''} \\
 & & & {\bullet} & & {\bullet} \ar[ll]^>{e}
}
$$

\noindent By following the labelling procedure, we obtain the following labelling of the vertices of $\Omega''$.
$$
\xymatrix @R=4pt @C=4pt {
{\Omega''} \\
 & & & {2} & & {1} \ar[ll] \\
\\
 & & & {} & & {} \ar[ll]^{p''}
}
$$

Having now completed the labelling of $\Omega''$, we add 3 to each of the labels of the vertices of $\Omega''$, and assign the labels thus obtained to the corresponding vertices in $\Omega'$.
$$
\xymatrix @R=4pt @C=4pt {
{\Omega'} & {} \ar[dddddddd]_{p'} & & {1} \ar[dd] \\
\\
 & & & {2} \ar[dd] \\
\\
 & & & {3} \ar[drr] \\
 & & & & & {5} \ar[dll] & & {4} \ar[ll] \\
 & & & {\bullet} \ar[uu] \ar[dd]_>{c} \\
\\
 & {} & & {\bullet}
}
$$

\noindent We label the secondary vertex corresponding to 3 relative to $p'$ with the label 6. Then, starting from 6, we proceed along $p'$ labelling subsequent vertices in increments of one. This gives us the following labelling of the vertices of $\Omega'$.
$$
\xymatrix @R=4pt @C=4pt {
{\Omega'} & & & {1} \ar[dd] \\
\\
 & & & {2} \ar[dd] \\
\\
 & & & {3} \ar[drr] \\
 & & & & & {5} \ar[dll] & & {4} \ar[ll] \\
 & & & {6} \ar[uu] \ar[dd] \\
\\
 & & & {7}
}
$$

\noindent We now add 2 to each of the labels of the vertices of $\Omega'$, and assign the labels thus obtained to the corresponding vertices in $\Omega$.
$$
\xymatrix @R=4pt @C=4pt {
{\Omega} & & & {3} \ar[dd] \\
\\
 & & & {4} \ar[dd] \\
\\
 & & & {5} \ar[drr] \\
 & & & & & {7} \ar[dll] & & {6} \ar[ll] \\
 & & & {8} \ar[uu] \ar[dd] \\
\\
 & & & {9} \ar[ddr] \\
\\
{1} \ar[rr] & & {2} \ar[uur] & & {\bullet} \ar[ll] \ar[rr]_>{b} & & {\bullet} \\
\\
{} \ar[rrrrrr]_p & & & & & & {}
}
$$

\noindent Finally, by labelling the secondary vertex corresponding to 2 relative to $p$ with the label 10, and then proceeding (from 10) along $p$, we complete the labelling of $\Omega$.
$$
\xymatrix @R=4pt @C=4pt {
{\Omega} & & & {3} \ar[dd] \\
\\
 & & & {4} \ar[dd] \\
\\
 & & & {5} \ar[drr] \\
 & & & & & {7} \ar[dll] & & {6} \ar[ll] \\
 & & & {8} \ar[uu] \ar[dd] \\
\\
 & & & {9} \ar[ddr] \\
\\
{1} \ar[rr] & & {2} \ar[uur] & & {10} \ar[ll] \ar[rr] & & {11}
}
$$
\eex


\section{Explicit construction of companion bases in Dynkin type $A_n$}
\label{type A construction}

From this point onwards, we suppose that the vertices of the quiver $\Gamma$ have been labelled in accordance with the labelling procedure outlined in Section~\ref{labelling section} above. By a minor abuse of notation, we will often refer to the vertices of $\Gamma$ by their labels as well as treating the labels as numerical values.

In this section, we construct a companion basis for $\Gamma$. This companion basis can be immediately written down simply by looking at the labels of the vertices of $\Gamma$. We start by observing some properties of the labelled quiver $\Gamma$.

During the procedure for labelling the vertices of $\Gamma$, a number of strings in $\Gamma$ are considered. These strings determine the labelling of the vertices of $\Gamma$, and we therefore refer to them as the \textit{labelling strings} for the given labelling of $\Gamma$. (The labelling strings are precisely the strings joining the pairs of end vertices considered during the labelling procedure.) The following result is clear.

\blemma
\label{properties of labelled quiver}
The labelled quiver $\Gamma$ has the following properties:\\
(i) Each vertex of $\Gamma$ lies on precisely one labelling string.\\
(ii) For any given 3-cycle in $\Gamma$, there is a unique labelling string which passes through exactly two vertices of that 3-cycle.\\
(iii) Suppose the vertex labelled $j$ is a primary vertex relative to some labelling string $p$, and that there are $b$ vertices in the subquiver of $\Gamma$ rooted at the complementary vertex corresponding to $j$ relative to $p$. Then, due to the inductive nature of the labelling procedure, we see that the complementary vertex corresponding to $j$ relative to $p$ will be labelled $j+b$. Also, by construction, the secondary vertex corresponding to $j$ relative to $p$  will be labelled $j+b+1$.
\elemma

In view of these properties, we make the following definition, in which reference to specific labelling strings for $\Gamma$ is dropped.

\bdfn
\label{3-cycle vertex types}
If the vertices of a 3-cycle in the labelled quiver $\Gamma$ have labels $i$, $j$ and $k$ with $i<j<k$, then we call $i$ a \textit{primary vertex}, $j$ a \textit{complementary vertex}, and $k$ a \textit{secondary vertex}.

Furthermore, we call $j$ the \textit{complementary vertex corresponding to $i$}, and we call $k$ the \textit{secondary vertex corresponding to $i$}.
\edfn

In the situation where two 3-cycles meet in a vertex, we see that that vertex can be both a secondary vertex and either a primary or a complementary vertex.

The following result is an immediate consequence of Lemma~\ref{properties of labelled quiver}.

\bcor
\label{left iff left rel path}
A given vertex of $\Gamma$ is a primary (resp.\ complementary, secondary) vertex if and only if it is a primary (resp.\ complementary, secondary) vertex relative to some labelling string for the given labelling of $\Gamma$.
\ecor

We are now able to state the result exhibiting a companion basis for $\Gamma$.

\bthm
\label{companion basis type A}
Let $\Gamma$ be a quiver of mutation type $A_n$, labelled according to the specified labelling procedure. The set $\{ \beta_i \colon i \in \Gamma_0 \} \subseteq \Phi$ given by
$$
\beta_i = \left \{
\begin{array}{cl}
\alpha_i + \ldots + \alpha_{i+a} & \textrm{if } i \textrm{ is a primary vertex and } i+a \textrm{ is the}\\
 & \textrm{complementary vertex corresponding to } i,\\
 & \\
\alpha_i & \textrm{if } i \textrm{ is not a primary vertex}
\end{array}
\right.
$$
is a companion basis for $\Gamma$.
\ethm

Note that for each $i$, $\beta_i$ is a root, since it is a sum of consecutive simple roots.

The remainder of this section is devoted to proving the result of Theorem~\ref{companion basis type A}. There are a couple of options available to us. Since the set specified in Theorem~\ref{companion basis type A} is clearly a $\mathbb{Z}$-basis for $\mathbb{Z} \Phi$, one option would be to further examine the properties of the labelled quiver $\Gamma$, and use these to directly check that the specified set really is a companion basis by evaluating inner products. This is the approach taken in \cite[Chapter 5]{Par1}. Here we take a different approach, calling upon the process of companion basis mutation introduced in \cite[Theorem 6.1]{Par2} (and \cite[Theorem 3.1.4]{Par1}). We feel that this is a more illuminating approach.

We start by noting that in order to find a companion basis for $\Gamma$, it is enough to find a companion basis for any quiver with the same underlying graph as $\Gamma$.

\bdfn
\label{quasi-isomorphism of quivers}
We say that two (labelled) quivers are \textit{quasi-isomorphic} if their underlying (labelled) graphs are isomorphic. In the case of labelled quivers, this means that their underlying graphs are isomorphic and corresponding vertices have the same label.
\edfn

The following is an immediate consequence of the definition of companion bases.

\blemma
\label{companion bases for quasi-isomorphic quivers}
Any companion basis for a quiver that is quasi-isomorphic to $\Gamma$ is also a companion basis for $\Gamma$ (considering the elements of such a companion basis to be indexed by the corresponding vertices of $\Gamma$).
\elemma

\brk
\label{companion bases for quasi-isomorphic quivers remark}
Despite the fact that we are restricting to the Dynkin type $A$ case in this section, we note that the result of Lemma~\ref{companion bases for quasi-isomorphic quivers} clearly holds in the more general case where we consider $\Gamma$ to be an arbitrary quiver of simply-laced Dynkin mutation type.
\erk

With the following result, we now show that by applying quiver mutations to the labelled Dynkin quiver
$$
\overrightarrow{A}_n \colon \; 1 \longrightarrow 2 \longrightarrow 3 \longrightarrow \cdots \longrightarrow n-1 \longrightarrow n,
$$
we are able to construct a labelled quiver quasi-isomorphic to $\Gamma$. In addition, the proof shows how the labelling of $\Gamma$ that we obtain using the prescribed labelling procedure, and which may seem a little unnatural at first sight, actually arises in a very natural way.

\bprop
\label{labelling by mutation}
We can obtain a labelled quiver quasi-isomorphic to $\Gamma$ from the Dynkin quiver $\overrightarrow{A}_n$ (labelled as shown above), by performing a sequence of quiver mutations, not requiring mutation at the vertices labelled 1 or $n$.
\eprop

\begin{proof}
We use induction on the number of vertices of $\Gamma$.

If $\Gamma$ has a single vertex (that is, if $n=1$), the result clearly holds.

Let $n \geq 2$ and suppose that the (corresponding) result holds for all labelled quivers of mutation type $A$ with fewer than $n$ vertices.

As we move along the string $p$ in $\Gamma$ from 1 to $n$, suppose that $i$ is the first primary vertex on $p$. Up to and including $i$, vertex labels increase in increments of 1. Note that there may be no primary vertices relative to $p$, in which case, the labelled quivers $\Gamma$ and $\overrightarrow{A}_n$ are quasi-isomorphic. Suppose that there are $a$ vertices in the subquiver of $\Gamma$ rooted at the complementary vertex corresponding to $i$ (so that this vertex must be labelled $i+a$). Mutating $\overrightarrow{A}_n$ successively at the vertices labelled $i+1 , \ldots , i+a$, we obtain the following:

$$
\xymatrix @R=4pt @C=4pt {
 & & & & & & & & & {i+1} \ar[dd] \\
\\
 & & & & & & & & & {i+2} \ar[dd] & & & & & & & & & {(\star)}\\
\\
 & & & & & & & & & {} \ar@{.}[dd] \\
\\
 & & & & & & & & & {} \ar[dd] \\
\\
 & & & & & & & & & {i+a} \ar[ddl] \\
\\
{1} \ar[rr] & & {2} \ar[rr] & & {} \ar@{.}[rr] & & {} \ar[rr] & & {i} \ar[rr] & & {i+a+1} \ar[rr] \ar[uul] & & {} \ar@{.}[rr] & & {} \ar[rr] & & {n-1} \ar[rr] & & {n} \\
}
$$

By construction, the full subquiver of this quiver on the vertices $1, \ldots, i$ and the full subquiver of $\Gamma$ on these vertices are quasi-isomorphic labelled quivers.

Consider the subquiver $\Gamma'$ of $\Gamma$ rooted at the complementary vertex $i+a$. The vertices of $\Gamma'$ are labelled with the labels $i+1, \ldots, i+a$. Due to the way in which the labelling of $\Gamma$ was produced, this labelling of $\Gamma'$ can be obtained by applying the labelling procedure to $\Gamma'$ and then adding $i$ to each vertex label. By induction, it therefore follows that we can obtain a labelled quiver quasi-isomorphic to $\Gamma'$ from the labelled quiver
$$
i+1 \longrightarrow \cdots \longrightarrow i+a
$$
by performing a sequence of quiver mutations (without any mutations at $i+1$ or $i+a$). Moreover, since mutation at $i+a$ is never required, these mutations can be performed within the labelled quiver in $(\star)$ without ever affecting the full subquiver on the vertices $\{ 1, \ldots, i, i+a+1, \ldots, n \}$.

Therefore, from the labelled Dynkin quiver $\overrightarrow{A}_n$, we are able to obtain a labelled quiver for which the labelled full subquiver on the vertices $1, \ldots, i+a+1$ is quasi-isomorphic to the labelled full subquiver of $\Gamma$ on the same vertices.

We can repeat this for each subsequent primary vertex lying on the string in $\Gamma$ from 1 to $n$.
\end{proof}

Denote the labelled quiver quasi-isomorphic to $\Gamma$ obtained by applying Proposition~\ref{labelling by mutation} to $\overrightarrow{A}_n$ by $\widetilde{\Gamma}$. The proof of Proposition~\ref{labelling by mutation} actually provides us with an explanation of how to construct the quiver $\widetilde{\Gamma}$, starting from the quiver $\overrightarrow{A}_n$ and applying quiver mutations. Since the simple system $\Pi = \{ \alpha_1 , \ldots , \alpha_n \}$ is a companion basis for $\overrightarrow{A}_n$ (where for each vertex $i$, $\alpha_i$ is the associated companion basis element), we are therefore able to obtain a companion basis for $\widetilde{\Gamma}$ by applying the corresponding sequence of companion basis mutations to $\Pi$. In order to do this, we must now recall the result of \cite{Par2} in which companion basis mutation was introduced.

Let $Q$ be a quiver of simply-laced Dynkin mutation type, and let $\Phi'$ be the associated root system. Suppose that $k$ is a vertex of $Q$, and that $Q'$ is the quiver obtained from $Q$ by applying quiver mutation at $k$. Denote the sets of vertices of $Q$ and $Q'$ by $Q_0$ and $Q_0'$ respectively. The following result \cite[Theorem 6.1]{Par2} shows how a companion basis for $Q$ can be mutated to produce a companion basis for $Q'$.

\bthm
\label{companion basis mutation theorem}
Let $\{ \gamma_x \colon x \in Q_0 \} \subseteq \Phi'$ be a companion basis for $Q$. Then,

(i) the set $\{ \gamma_x' \colon x \in Q_0' \} \subseteq \Phi'$ given by
$$
\gamma_x' = \left \{
\begin{array}{cl}
s_{\gamma_k} (\gamma_x) & \textrm{if there is an arrow from } x \textrm{ to } k \textrm{ in } Q,\\
\gamma_x & \textrm{otherwise}
\end{array}
\right.
$$
is a companion basis for $Q'$;

(ii) the set $\{ \gamma_x'' \colon x \in Q_0' \} \subseteq \Phi'$ given by
$$
\gamma_x'' = \left \{
\begin{array}{cl}
s_{\gamma_k} (\gamma_x) & \textrm{if there is an arrow from } k \textrm{ to } x \textrm{ in } Q,\\
\gamma_x & \textrm{otherwise}
\end{array}
\right.
$$
is a companion basis for $Q'$.

We refer to $\{ \gamma_x' \colon x \in Q_0' \}$ as the companion basis for $Q'$ obtained from $\{ \gamma_x \colon x \in Q_0 \}$ by mutating inwardly at $k$, and we refer to $\{ \gamma_x'' \colon x \in Q_0' \}$ as the companion basis for $Q'$ obtained from $\{ \gamma_x \colon x \in Q_0 \}$ by mutating outwardly at $k$.
\ethm

\brk
\label{comp basis mutation remark}
A benefit of having two types of companion basis mutation is that we can choose the one which makes life easiest for us in any given situation.

Barot and Marsh established an extension of Theorem~\ref{companion basis mutation theorem} to all finite type cases in \cite[Proposition 6.4]{BM}.
\erk

With the following result, we construct a companion basis for $\widetilde{\Gamma}$. This companion basis coincides with our candidate companion basis from the statement of Theorem~\ref{companion basis type A}. Therefore, since $\Gamma$ and $\widetilde{\Gamma}$ are quasi-isomorphic labelled quivers, this completes the proof of Theorem~\ref{companion basis type A}. Denote the set of vertices of $\widetilde{\Gamma}$ by $\widetilde{\Gamma}_0$.

\bprop
\label{companion basis for quasi-isomorphic quiver}
The set $\{ \beta_t \colon t \in \widetilde{\Gamma}_0 \} \subseteq \Phi$ given by
$$
\beta_t = \left \{
\begin{array}{cl}
\alpha_t + \ldots + \alpha_{t+a} & \textrm{if } t \textrm{ is a primary vertex and } t+a \textrm{ is the}\\
 & \textrm{complementary vertex corresponding to } t,\\
 & \\
\alpha_t & \textrm{if } t \textrm{ is not a primary vertex}
\end{array}
\right.
$$
is a companion basis for $\widetilde{\Gamma}$.
\eprop

\begin{proof}
We consider the procedure of constructing $\widetilde{\Gamma}$ from $\overrightarrow{A}_n$ by quiver mutations and examine the effect of performing the corresponding inward companion basis mutations, starting from the companion basis $\Pi$ for $\overrightarrow{A}_n$.

The first step in constructing the quiver $\widetilde{\Gamma}$ is to mutate $\overrightarrow{A}_n$ consecutively at the vertices $i+1 , \ldots , i+a$ where:\\
(i) $i$ is the first primary vertex we encounter when moving along the string in $\Gamma$ from 1 to $n$;\\
(ii) $i+a$ is the complementary vertex corresponding to $i$ in $\Gamma$.

Let $\widetilde{\Gamma}^{(j)}$ be the (labelled) quiver obtained from $\overrightarrow{A}_n$ by mutating consecutively at the vertices $i+1 , \ldots , i+j$. For convenience, we shall additionally use $\widetilde{\Gamma}^{(0)}$ to denote the (labelled) quiver $\overrightarrow{A}_n$.

Let $\{ \beta_t^{(j)} \colon t \in \widetilde{\Gamma}_0^{(j)} \} \subseteq \Phi$ be the companion basis for $\widetilde{\Gamma}^{(j)}$ obtained by applying inward companion basis mutation consecutively at the vertices $i+1 , \ldots , i+j$ to the companion basis $\Pi$ for $\overrightarrow{A}_n$. Note here that $\Pi = \{ \beta_t^{(0)} \colon t \in \widetilde{\Gamma}_0^{(0)} \}$, so that $\beta_t^{(0)} = \alpha_t$ for all $t \in \widetilde{\Gamma}_0^{(0)}$.

We now show that
$$
\beta_t^{(j)} = \left \{
\begin{array}{cl}
\alpha_i + \ldots + \alpha_{i+j} & \textrm{if } t = i\\
\alpha_t & \textrm{if } t \neq i
\end{array}
\right.
$$
by induction on $j$.

The result holds in the initial case where $j=0$.

Suppose the result also holds for $j-1$ and consider the quiver $\widetilde{\Gamma}^{(j-1)}$:
$$
\xymatrix @R=4pt @C=4pt {
{\widetilde{\Gamma}^{(j-1)}} & & & & & & & {i+1} \ar[dd] \\
\\
 & & & & & & & {} \ar@{.}[dd] \\
\\
 & & & & & & & {} \ar[dd] \\
\\
 & & & & & & & {i+j-1} \ar[ddl] \\
\\
{1} \ar[rr] & & {} \ar@{.}[rr] & & {} \ar[rr] & & {i} \ar[rr] & & {i+j} \ar[rr] \ar[uul] & & {} \ar@{.}[rr] & & {} \ar[rr] & & {n} \\
}
$$

We get the quiver $\widetilde{\Gamma}^{(j)}$ by mutating $\widetilde{\Gamma}^{(j-1)}$ at the vertex $i+j$. And we get our desired companion basis for $\widetilde{\Gamma}^{(j)}$ by mutating the companion basis (provided by the induction hypothesis) for $\widetilde{\Gamma}^{(j-1)}$ inwardly at $i+j$. The only arrow with head $i+j$ in $\widetilde{\Gamma}^{(j-1)}$ has source $i$, so applying companion basis mutation, we get:
\begin{eqnarray*}
\beta_t^{(j)} & = & \beta_t^{(j-1)} = \alpha_t, \quad t \neq i \\
\beta_i^{(j)} & = & s_{\beta_{i+j}^{(j-1)}} (\beta_i^{(j-1)}) = s_{\alpha_{i+j}} (\alpha_i + \cdots + \alpha_{i+j-1}) = \alpha_i + \cdots + \alpha_{i+j}
\end{eqnarray*}
as required. (Note that this also works when passing from $\widetilde{\Gamma}^{(0)}$ to $\widetilde{\Gamma}^{(1)}$.) In particular, the companion basis we obtain for $\widetilde{\Gamma}^{(a)}$ by mutating inwardly at $i+1 , \ldots , i+a$ is
$$
\beta_t^{(a)} = \left \{
\begin{array}{cl}
\alpha_i + \ldots + \alpha_{i+a} & \textrm{if } t = i\\
\alpha_t & \textrm{if } t \neq i.
\end{array}
\right.
$$

$$
\xymatrix @R=4pt @C=4pt {
{\widetilde{\Gamma}^{(a)}} & & & & & & & {i+1} \ar[dd] \\
\\
 & & & & & & & {} \ar@{.}[dd] \\
\\
 & & & & & & & {} \ar[dd] \\
\\
 & & & & & & & {i+a} \ar[ddl] \\
\\
{1} \ar[rr] & & {} \ar@{.}[rr] & & {} \ar[rr] & & {i} \ar[rr] & & {i+a+1} \ar[rr] \ar[uul] & & {} \ar@{.}[rr] & & {} \ar[rr] & & {n} \\
}
$$

The result follows due to the inductive nature of the procedure in which $\widetilde{\Gamma}$ is constructed from $\overrightarrow{A}_n$. (It's important to note here that for each vertex of $\widetilde{\Gamma}^{(a)}$ in $\{ i+1, \ldots, n \}$ the associated root in the companion basis we obtained for $\widetilde{\Gamma}^{(a)}$ is the corresponding simple root.)
\end{proof}

We conclude this section with an example which, when considered together with Example~\ref{labelling example}, serves to illustrate the ease with which our companion basis construction procedure can be applied in practice.

\bex
\label{explicit beta set example}
Using Theorem~\ref{companion basis type A}, we write down a companion basis for the labelled quiver $\Omega$ of mutation type $A_{11}$ considered in Example~\ref{labelling example}:
$$
\xymatrix @R=4pt @C=4pt {
{\Omega} & & & {3} \ar[dd] \\
\\
 & & & {4} \ar[dd] \\
\\
 & & & {5} \ar[drr] \\
 & & & & & {7} \ar[dll] & & {6} \ar[ll] \\
 & & & {8} \ar[uu] \ar[dd] \\
\\
 & & & {9} \ar[ddr] \\
\\
{1} \ar[rr] & & {2} \ar[uur] & & {10} \ar[ll] \ar[rr] & & {11}
}
$$

Take $\Phi$ to be the root system of Dynkin type $A_{11}$, with simple system $\Pi = \{ \alpha_1 , \ldots , \alpha_{11} \}$.

We start by noting that 2 and 5 are the only primary vertices in $\Omega$. We therefore set $\beta_i = \alpha_i$ for $i \neq 2,5$. The complementary vertex corresponding to 2 is 9, so we set $\beta_2 = \alpha_2 + \cdots + \alpha_9$. Likewise, the complementary vertex corresponding to 5 is 7 and so we set $\beta_5 = \alpha_5 + \alpha_6 + \alpha_7$. We then have that $\{ \beta_1 , \ldots , \beta_{11} \} \subseteq \Phi$ is a companion basis for $\Omega$ (we note that this can be easily verified directly).
\eex


\section{Explicit construction of companion bases in Dynkin type $D_n$}
\label{type D construction}

We now turn and restrict our attention to the case where $\Gamma$ is a quiver of mutation type $D_n$, $n\geq 4$. We show how to construct an explicit companion basis for $\Gamma$. The quivers of mutation type $D_n$ were described by Vatne in \cite{Vat} and are in fact closely related to the quivers of mutation type $A$. We are able to take advantage of this fact in our companion basis construction procedure for type $D$, which makes heavy use of that for type $A$. The operation of companion basis mutation again plays an important role. As in the type $A$ case, the constructed companion basis for $\Gamma$ can simply be written down upon looking at the structure of $\Gamma$ (that is, it can be read-off from the quiver, after applying a labelling).

We start by recalling Vatne's description of the quivers of mutation type $D_n$. This description tells us that the quiver $\Gamma$ is of one of four possible types. Each type consists of a collection of quivers of mutation type $A$ `glued' to a `skeleton' quiver (using the terminology of \cite[Section 1.6]{BHL}). The four types are outlined below.

\underline{Type I}: Let $\Gamma^{(0)}$ be the following quiver

 \begin{center}
 \psfragscanon
 \psfrag{a}{\mbox{\small $a$}}
 \psfrag{b}{\mbox{\small $b$}}
 \psfrag{c}{\mbox{\small $c_1$}}
 \psfrag{g0}{$\Gamma^{(0)}$:}
 \includegraphics[scale=.50]{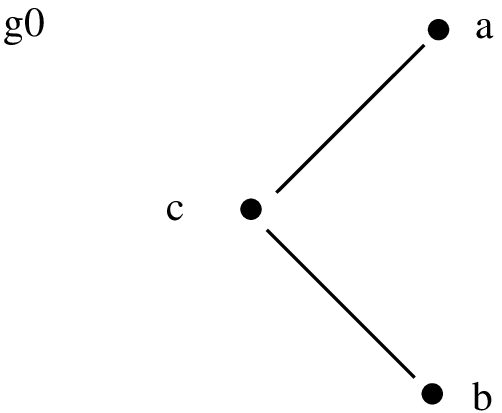}
 \end{center}

\noindent where the edges joining $c_1$ to $a,b$ can have any orientation. In this case, $\Gamma$ is obtained from a disjoint union $\Gamma^{(0)} \sqcup \Gamma^{(1)}$, where $\Gamma^{(1)}$ is a quiver of mutation type $A_{n-2}$, by identifying some end vertex $v_1$ of $\Gamma^{(1)}$ with the vertex $c_1$. This is illustrated in the following diagram:

 \begin{center}
 \psfragscanon
 \psfrag{a}{\mbox{\small $a$}}
 \psfrag{b}{\mbox{\small $b$}}
 \psfrag{c}{\mbox{\small $v_1 \equiv c_1$}}
 \psfrag{g1}{\mbox{\small $\Gamma^{(1)}$}}
 \psfrag{g}{$\Gamma$:}
 \includegraphics[scale=.50]{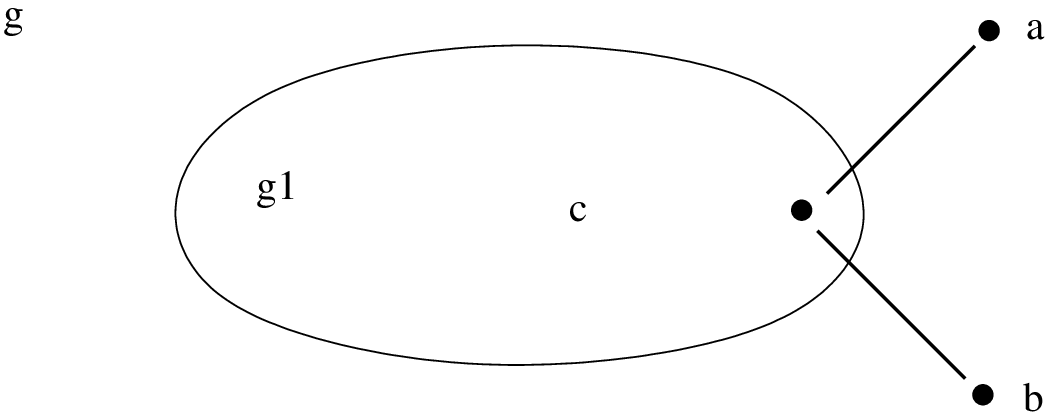}
 \end{center}

\underline{Type II}: Let $\Gamma^{(0)}$ be the following quiver:

 \begin{center}
 \psfragscanon
 \psfrag{c1}{\mbox{\small $c_1$}}
 \psfrag{c2}{\mbox{\small $c_2$}}
 \psfrag{g0}{$\Gamma^{(0)}$:}
 \includegraphics[scale=.50]{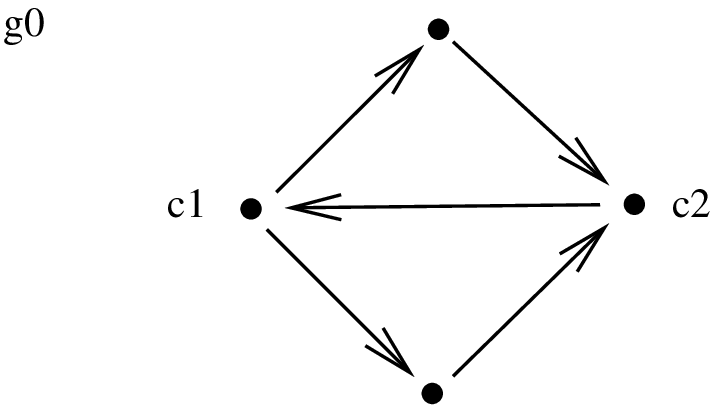}
 \end{center}

In this case, the quiver $\Gamma$ is obtained from a disjoint union $\Gamma^{(0)} \sqcup \Gamma^{(1)} \sqcup \Gamma^{(2)}$, where $\Gamma^{(1)}$ and $\Gamma^{(2)}$ are quivers of mutation types $A_m$ and $A_{n-m-2}$ respectively, for some $1 \leq m \leq n-3$. Indeed, $\Gamma$ is obtained from this disjoint union by identifying some end vertices $v_1$ of $\Gamma^{(1)}$ and $v_2$ of $\Gamma^{(2)}$ with the vertices $c_1$ and $c_2$ respectively. This is illustrated in the following diagram:

 \begin{center}
 \psfragscanon
 \psfrag{c1}{\mbox{\small $v_1 \equiv c_1$}}
 \psfrag{c2}{\mbox{\small $v_2 \equiv c_2$}}
 \psfrag{g1}{\mbox{\small $\Gamma^{(1)}$}}
 \psfrag{g2}{\mbox{\small $\Gamma^{(2)}$}}
 \psfrag{g}{$\Gamma$:}
 \includegraphics[scale=.50]{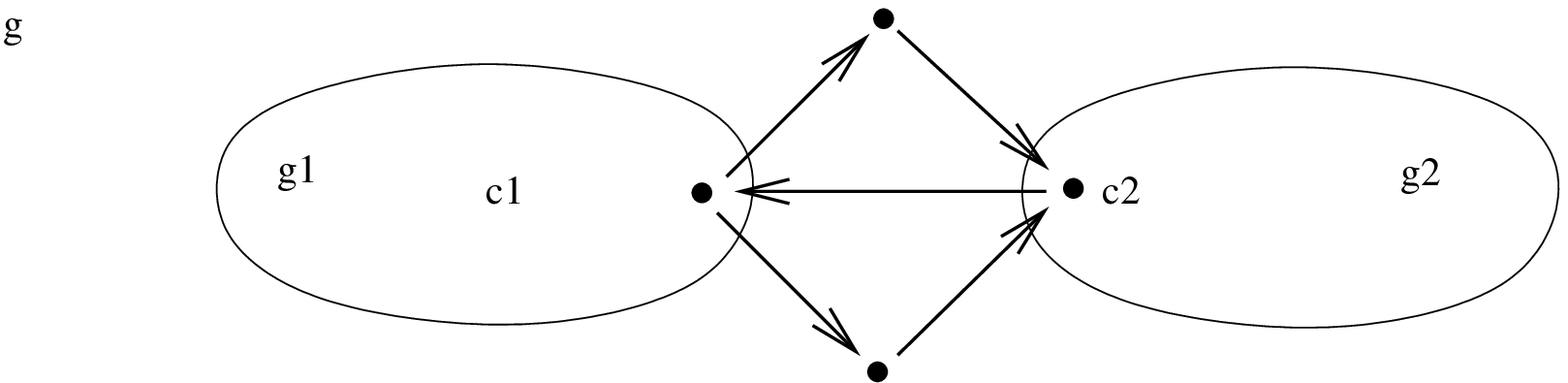}
 \end{center}

\underline{Type III}: Let $\Gamma^{(0)}$ be the following quiver:

 \begin{center}
 \psfragscanon
 \psfrag{c1}{\mbox{\small $c_1$}}
 \psfrag{c2}{\mbox{\small $c_2$}}
 \psfrag{g0}{$\Gamma^{(0)}$:}
 \includegraphics[scale=.50]{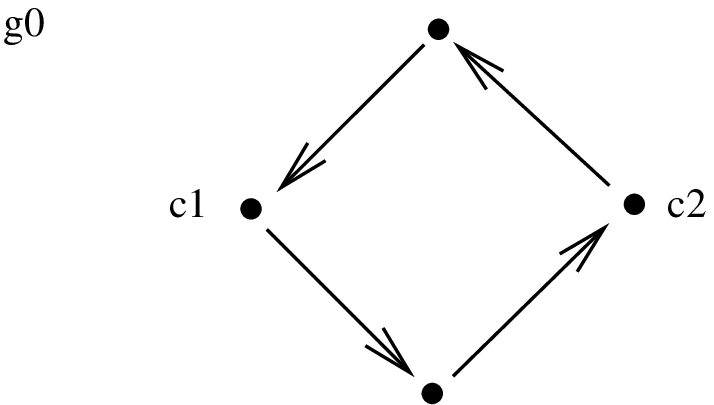}
 \end{center}

In this case, the quiver $\Gamma$ is obtained from a disjoint union $\Gamma^{(0)} \sqcup \Gamma^{(1)} \sqcup \Gamma^{(2)}$, where $\Gamma^{(1)}$ and $\Gamma^{(2)}$ are quivers of mutation types $A_m$ and $A_{n-m-2}$ respectively, for some $1 \leq m \leq n-3$. Indeed, $\Gamma$ is obtained from this disjoint union by identifying some end vertices $v_1$ of $\Gamma^{(1)}$ and $v_2$ of $\Gamma^{(2)}$ with the vertices $c_1$ and $c_2$ respectively. This is illustrated in the following diagram:

 \begin{center}
 \psfragscanon
 \psfrag{c1}{\mbox{\small $v_1 \equiv c_1$}}
 \psfrag{c2}{\mbox{\small $v_2 \equiv c_2$}}
 \psfrag{g1}{\mbox{\small $\Gamma^{(1)}$}}
 \psfrag{g2}{\mbox{\small $\Gamma^{(2)}$}}
 \psfrag{g}{$\Gamma$:}
 \includegraphics[scale=.50]{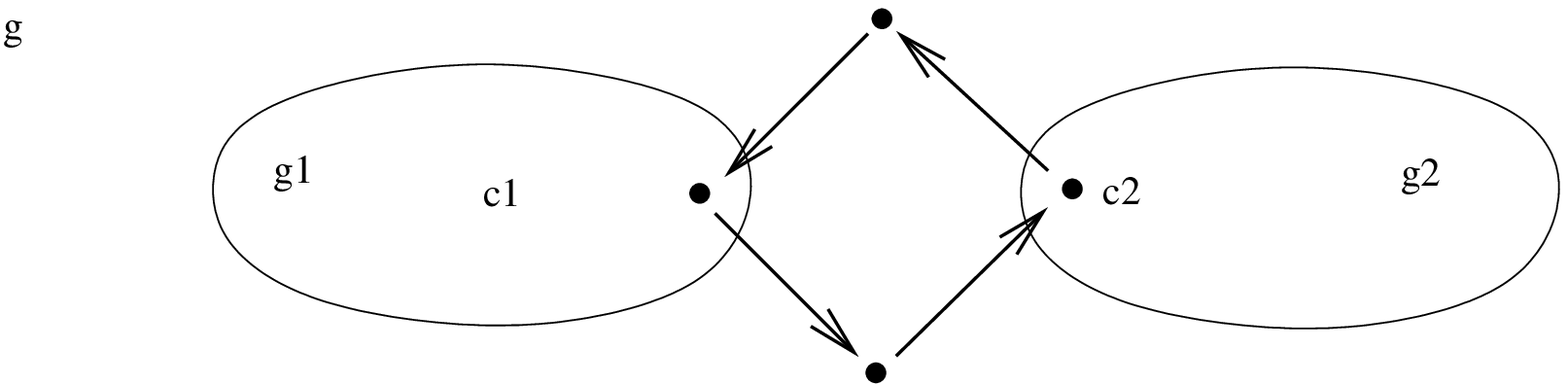}
 \end{center}

\underline{Type IV}: In this case, $\Gamma$ has a full subquiver $\Gamma^{(0)}$ (the skeleton quiver here), which we start by describing. $\Gamma^{(0)}$ contains a \textit{central cycle}, which is an oriented $m$-cycle, for some $m \geq 3$. Moreover, $\Gamma^{(0)}$ has $r$, $0 \leq r \leq m$, additional vertices $c_1, \ldots c_r$, which correspond respectively to arrows $p_1 \stackrel{\alpha_1}{\longrightarrow} s_1 , \ldots ,p_r \stackrel{\alpha_r}{\longrightarrow} s_r$ appearing in clockwise order around the central cycle. For each such pair $\alpha_i$ and $c_i$, there is an oriented 3-cycle $p_i \stackrel{\alpha_i}{\longrightarrow} s_i \longrightarrow c_i \longrightarrow p_i$, which we call a \textit{spike}. Any number of the arrows on the central cycle may belong to a spike. However, if $m=3$, we must have $r>0$. An example of a type IV skeleton quiver is given in the following diagram:

 \begin{center}
 \psfragscanon
 \psfrag{p1}{\mbox{\small $p_1$}}
 \psfrag{s1}{\mbox{\small $s_1$}}
 \psfrag{c1}{\mbox{\small $c_1$}}
 \psfrag{a1}{\mbox{\small $\alpha_1$}}
 \psfrag{p2}{\mbox{\small $p_2$}}
 \psfrag{s2p3}{\mbox{\small $s_2 = p_3$}}
 \psfrag{c2}{\mbox{\small $c_2$}}
 \psfrag{a2}{\mbox{\small $\alpha_2$}}
 \psfrag{s3}{\mbox{\small $s_3$}}
 \psfrag{c3}{\mbox{\small $c_3$}}
 \psfrag{a3}{\mbox{\small $\alpha_3$}}
 \psfrag{pr}{\mbox{\small $p_r$}}
 \psfrag{sr}{\mbox{\small $s_r$}}
 \psfrag{cr}{\mbox{\small $c_r$}}
 \psfrag{ar}{\mbox{\small $\alpha_r$}}
 \psfrag{G0}{$\Gamma^{(0)}$:}
 \includegraphics[scale=.50]{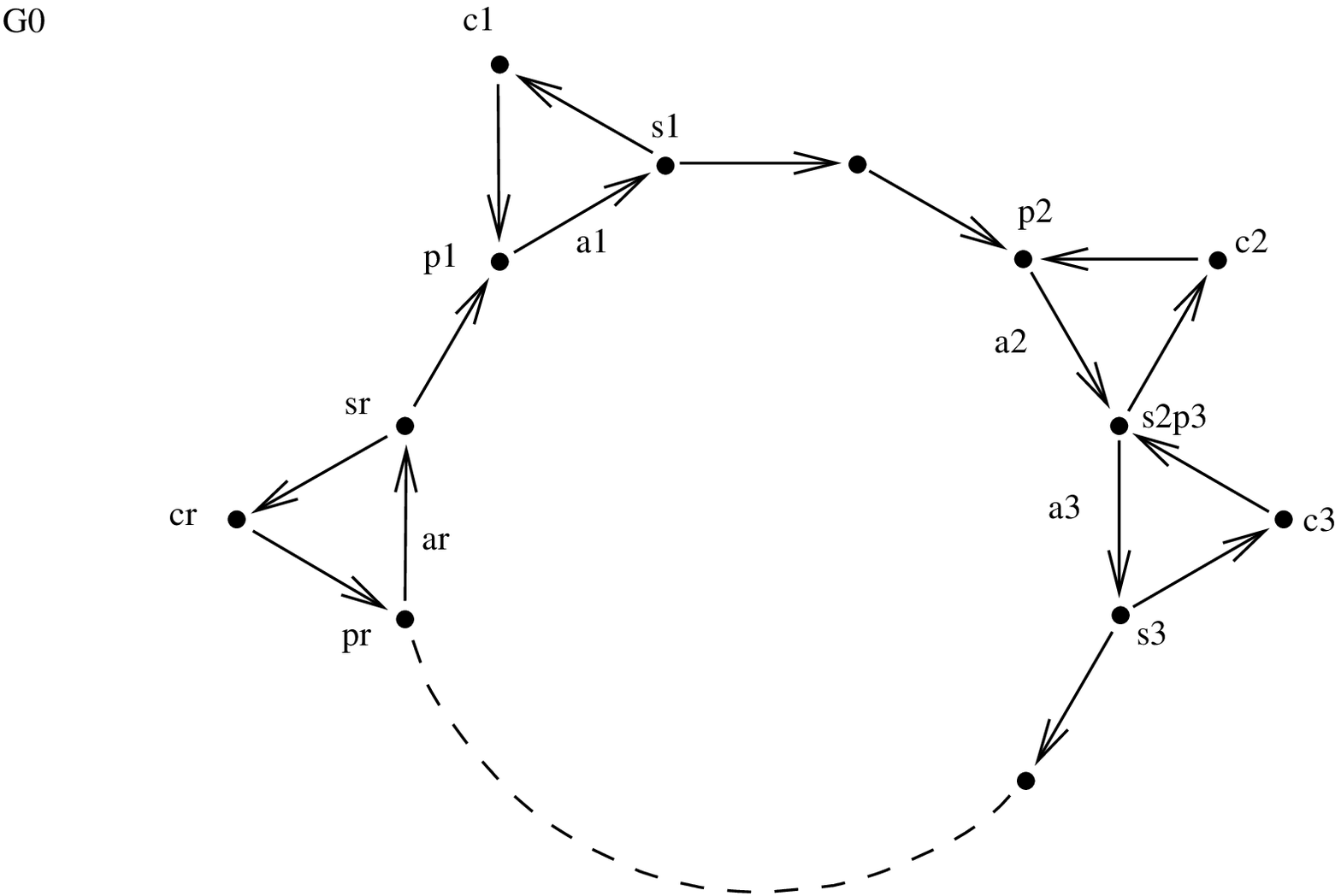}
 \end{center}

\noindent $\Gamma$ is then obtained from a disjoint union $\Gamma^{(0)} \sqcup \Gamma^{(1)} \sqcup \cdots \sqcup \Gamma^{(r)}$, for some collection $\Gamma^{(1)}, \ldots , \Gamma^{(r)}$ of quivers of mutation type $A$, by identifying some collection $v_1, \ldots , v_r$ of end vertices of $\Gamma^{(1)}, \ldots , \Gamma^{(r)}$ respectively with the vertices $c_1 , \ldots , c_r$ respectively. This is illustrated in the following diagram:

 \begin{center}
 \psfragscanon
 \psfrag{p1}{\mbox{\small $p_1$}}
 \psfrag{s1}{\mbox{\small $s_1$}}
 \psfrag{v1c1}{\mbox{\small $v_1 \equiv c_1$}}
 \psfrag{a1}{\mbox{\small $\alpha_1$}}
 \psfrag{p2}{\mbox{\small $p_2$}}
 \psfrag{s2p3}{\mbox{\small $s_2 = p_3$}}
 \psfrag{v2c2}{\mbox{\small $v_2 \equiv c_2$}}
 \psfrag{a2}{\mbox{\small $\alpha_2$}}
 \psfrag{s3}{\mbox{\small $s_3$}}
 \psfrag{v3c3}{\mbox{\small $v_3 \equiv c_3$}}
 \psfrag{a3}{\mbox{\small $\alpha_3$}}
 \psfrag{pr}{\mbox{\small $p_r$}}
 \psfrag{sr}{\mbox{\small $s_r$}}
 \psfrag{vrcr}{\mbox{\small $v_r \equiv c_r$}}
 \psfrag{ar}{\mbox{\small $\alpha_r$}}
 \psfrag{G}{$\Gamma$:}
 \psfrag{G1}{\mbox{\small $\Gamma^{(1)}$}}
 \psfrag{G2}{\mbox{\small $\Gamma^{(2)}$}}
 \psfrag{G3}{\mbox{\small $\Gamma^{(3)}$}}
 \psfrag{Gr}{\mbox{\small $\Gamma^{(r)}$}}
 \includegraphics[scale=.50]{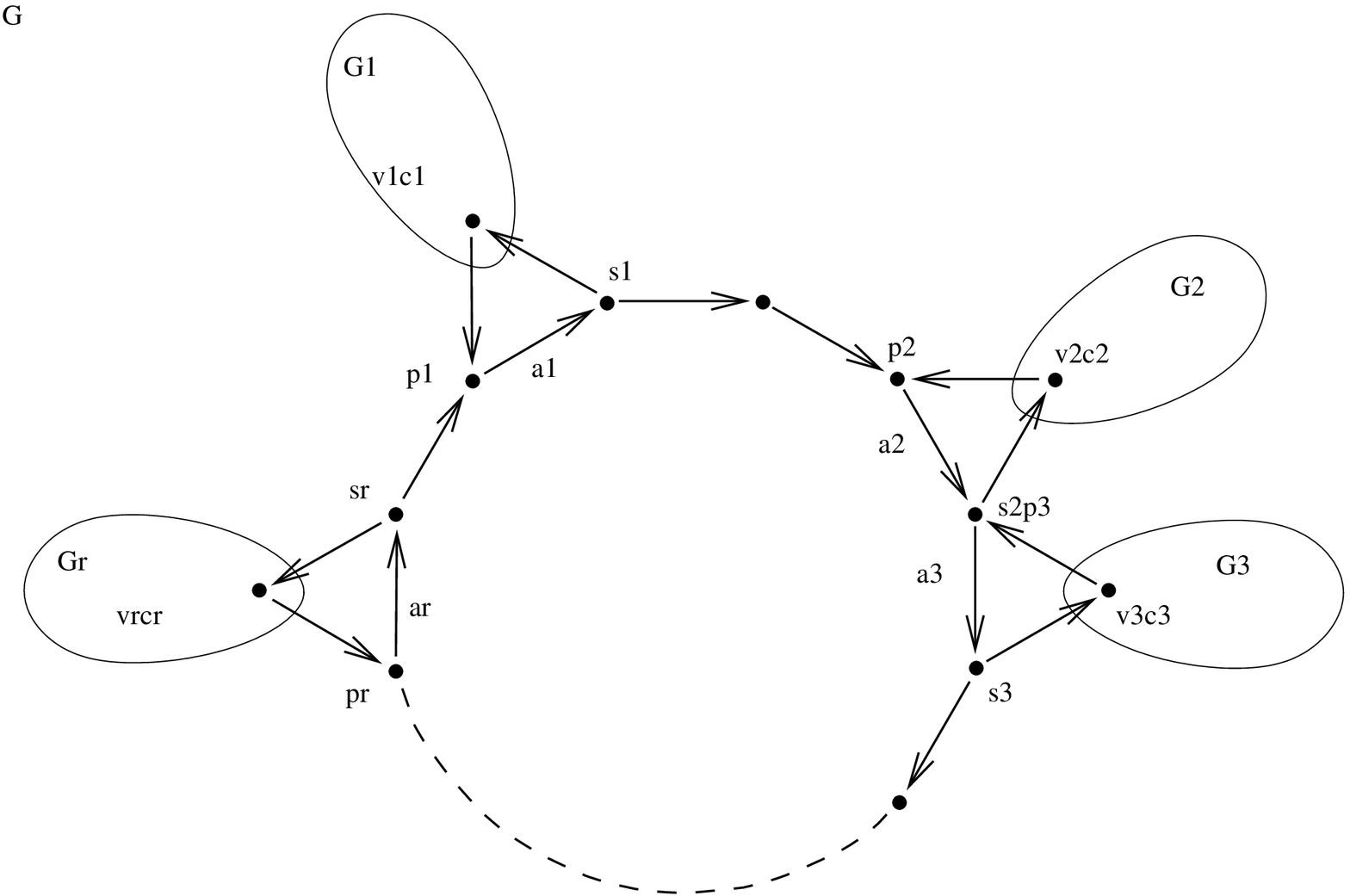}
 \end{center}

So, in each type, there is some quiver $\Gamma^{(0)}$, some $r \geq 0$, and some collection $\Gamma^{(1)}, \ldots , \Gamma^{(r)}$ of quivers of mutation type $A$ (this collection is empty when $r=0$) such that $\Gamma$ is obtained from a disjoint union $\Gamma^{(0)} \sqcup \Gamma^{(1)} \sqcup \cdots \sqcup \Gamma^{(r)}$ by identifying an end vertex $v_i$ of $\Gamma^{(i)}$ with a distinct vertex $c_i$ of $\Gamma^{(0)}$ for each $i$, $1 \leq i \leq r$. We write $\Gamma = \Gamma^{(0)} \underline{\sqcup} \left( \Gamma^{(1)}, \ldots, \Gamma^{(r)} \right)$ to denote this ($r \neq 0$).

We consider each type separately, constructing a companion basis for $\Gamma$ in each case. In each case, we first label the vertices of $\Gamma$ and explain how a labelled quiver quasi-isomorphic to $\Gamma$ can be obtained from the labelled quiver

 \begin{center}
 \psfragscanon
 \psfrag{1}{\mbox{\small $1$}}
 \psfrag{2}{\mbox{\small $2$}}
 \psfrag{n-3}{\mbox{\small $n-3$}}
 \psfrag{n-2}{\mbox{\small $n-2$}}
 \psfrag{n-1}{\mbox{\small $n-1$}}
 \psfrag{n}{\mbox{\small $n$}}
 \psfrag{dn}{$\overrightarrow{D}_n$:}
 \includegraphics[scale=.50]{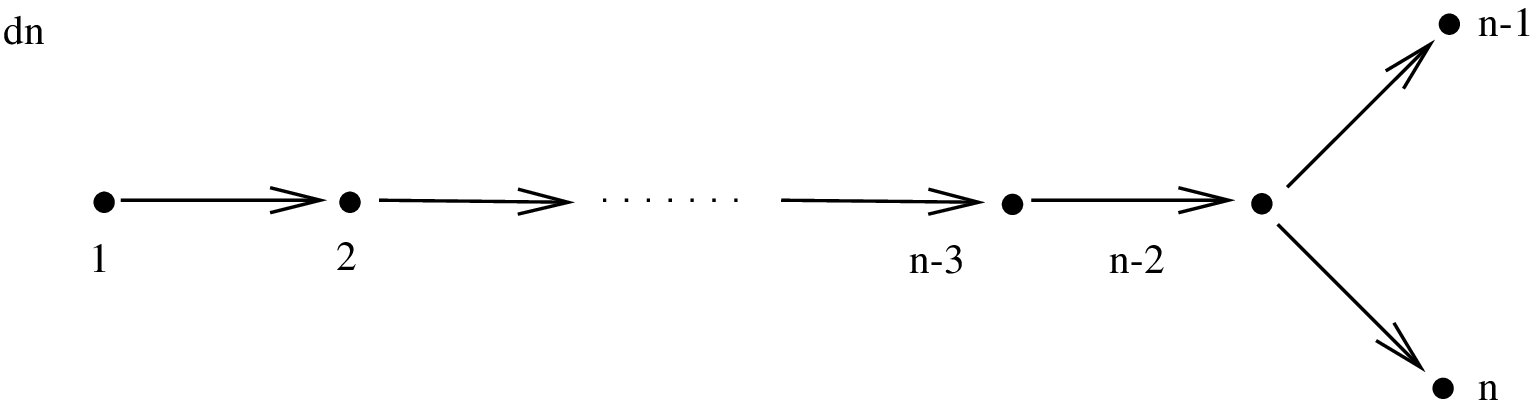}
 \end{center}

\noindent by applying a sequence of quiver mutations. Starting from the companion basis $\Pi = \{ \alpha_1 , \ldots , \alpha_n \}$ for $\overrightarrow{D}_n$, we then construct our explicit companion basis for $\Gamma$ by applying a corresponding sequence of companion basis mutations. For the full subquivers $\Gamma^{(1)}, \ldots , \Gamma^{(r)}$ of $\Gamma$, we can use our type $A$ companion basis construction procedure (which hides some of the details). We are able to do this due to the fact that the root system of type $A_m$, for $m<n$, naturally embeds into the root system of type $D_n$.

With the following definition, we extend the notion of primary, secondary and complementary vertices to the type $D$ case. These will again play a central role in enabling us to read-off a companion basis for $\Gamma$, once a labelling has been applied.

\bdfn
\label{vertex types in type D}
In the labelled quiver $\Gamma$, we call a vertex a \textit{primary} (resp. \textit{secondary}, \textit{complementary}) \textit{vertex} if it is a primary (resp. secondary, complementary) vertex in $\Gamma^{(i)}$ for some $i$, $1 \leq i \leq r$.
\edfn

\underline{Type I}: We have that $\Gamma = \Gamma^{(0)} \underline{\sqcup} \left ( \Gamma^{(1)} \right )$, with the end vertex $v_1$ of the quiver $\Gamma^{(1)}$ of mutation type $A_{n-2}$ being identified with the vertex $c_1$ of the skeleton quiver $\Gamma^{(0)}$, as shown previously. Now, we choose an ordered pair of distinct end vertices in $\Gamma^{(1)}$ with $v_1$ being the second vertex in this pair. Applying our labelling procedure from the previous section, we then obtain a labelling of $\Gamma^{(1)}$. We then apply the same labels to the corresponding vertices of $\Gamma$ and label the remaining two vertices $n-1$ and $n$, as shown below. (Note that the vertex $v_1$ gets the label $n-2$.)

 \begin{center}
 \psfragscanon
 \psfrag{n-2}{\mbox{\small $n-2$}}
 \psfrag{n-1}{\mbox{\small $n-1$}}
 \psfrag{n}{\mbox{\small $n$}}
 \psfrag{g1}{\mbox{\small $\Gamma^{(1)}$}}
 \psfrag{g}{$\Gamma$:}
 \includegraphics[scale=.50]{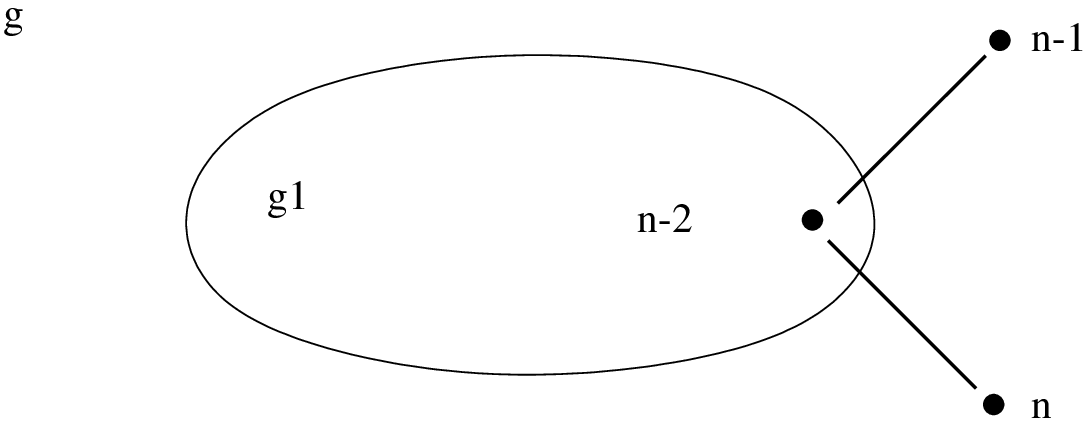}
 \end{center}

We saw in the previous section that by applying quiver mutations to the labelled quiver $\overrightarrow{A}_{n-2}$, we can obtain a labelled quiver $\widetilde{\Gamma}^{(1)}$ quasi-isomorphic to $\Gamma^{(1)}$. Moreover, since mutation at $n-2$ is not required, performing the corresponding sequence of mutations at the vertices of the labelled quiver $\overrightarrow{D}_n$, we obtain a labelled quiver $\widetilde{\Gamma}$ quasi-isomorphic to $\Gamma$. (Note that here we are identifying the quiver $\overrightarrow{A}_{n-2}$ with the full subquiver of $\overrightarrow{D}_n$ on the vertices $\{1, \ldots , n-2\}$.)

\bthm
\label{type D type I}
Let $\Gamma$ be a type I quiver of mutation type $D_n$, labelled as specified above. The set $\{ \beta_i \colon i \in \Gamma_0 \} \subseteq \Phi$ given by
$$
\beta_i = \left \{
\begin{array}{cl}
\alpha_i + \ldots + \alpha_{i+a} & \textrm{if } i \textrm{ is a primary vertex and } i+a \textrm{ is the}\\
 & \textrm{complementary vertex corresponding to } i,\\
 & \\
\alpha_i & \textrm{if } i \textrm{ is not a primary vertex}
\end{array}
\right.
$$
is a companion basis for $\Gamma$.
\ethm

\begin{proof}
Consider the sequence of quiver mutations taking us from $\overrightarrow{D}_n$ to $\widetilde{\Gamma}$. Starting from the companion basis $\Pi = \{ \alpha_1 , \ldots , \alpha_n \}$ for $\overrightarrow{D}_n$, we may obtain a companion basis for $\widetilde{\Gamma}$ (and hence also for $\Gamma$) by performing the corresponding sequence of inward companion basis mutations.

Since the root system of type $A_{n-2}$ naturally embeds into the root system of type $D_n$, we can consider the root system of type $A_{n-2}$ as a subset of the root system of type $D_n$. We may therefore proceed exactly as in the proof of Proposition~\ref{companion basis for quasi-isomorphic quiver} to see that the companion basis elements we obtain associated to the vertices of $\widetilde{\Gamma}$ belonging to $\widetilde{\Gamma}^{(1)}$ match those given in the statement of the theorem. Moreover, since no companion basis mutation is performed at the vertex $n-2$ (because no quiver mutation is), it follows immediately from Theorem~\ref{companion basis mutation theorem} that the companion basis elements we obtain associated to the vertices $n-1$ and $n$ are $\alpha_{n-1}$ and $\alpha_n$ respectively (as there can never be an arrow going from $n-1$ or $n$ into a vertex of mutation).
\end{proof}

\underline{Type II}: We have that $\Gamma = \Gamma^{(0)} \underline{\sqcup} \left( \Gamma^{(1)}, \Gamma^{(2)} \right)$, with the end vertex $v_1$ of $\Gamma^{(1)}$ and the end vertex $v_2$ of $\Gamma^{(2)}$ identified respectively with the vertices $c_1$ and $c_2$ of the skeleton quiver $\Gamma^{(0)}$, as shown previously. Here, $\Gamma^{(1)}$ is a quiver of mutation type $A_m$ and $\Gamma^{(2)}$ is a quiver of mutation type $A_{n-m-2}$ (for some $1 \leq m \leq n-3$).

Starting from $\overrightarrow{D}_n$ and applying consecutive quiver mutations at the vertices $n-2, n-3, \ldots , m+1$, we obtain the labelled quiver, which we call $\Xi_m$, as shown below:

 \begin{center}
 \psfragscanon
 \psfrag{1}{\mbox{\small $1$}}
 \psfrag{m}{\mbox{\small $m$}}
 \psfrag{m+1}{\mbox{\small $m+1$}}
 \psfrag{n-2}{\mbox{\small $n-2$}}
 \psfrag{n-1}{\mbox{\small $n-1$}}
 \psfrag{n}{\mbox{\small $n$}}
 \psfrag{d}{$\Xi_m$:}
 \includegraphics[scale=.50]{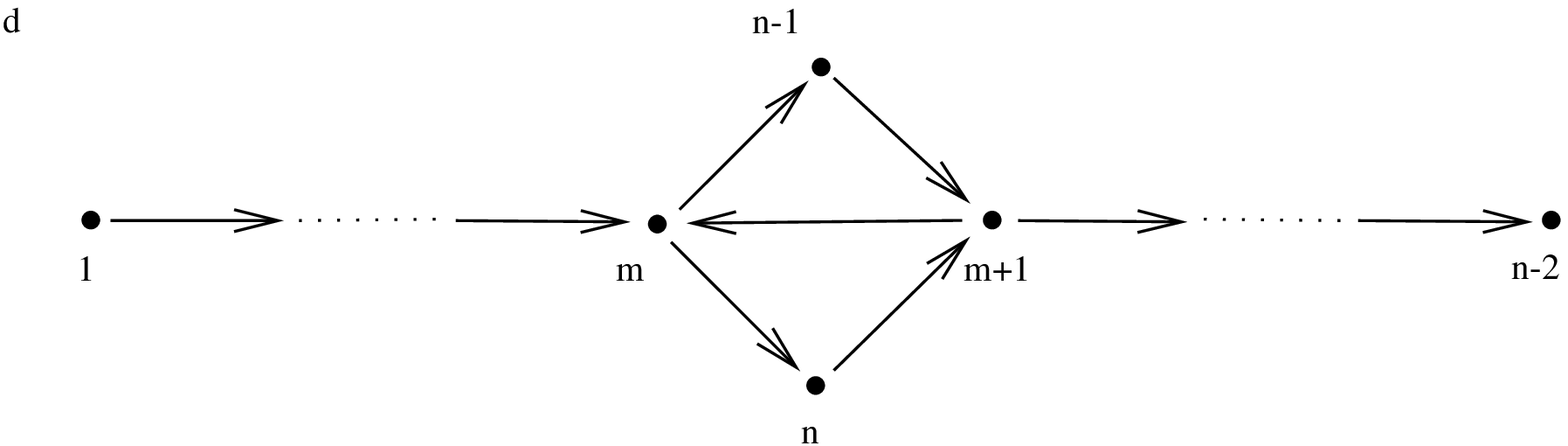}
 \end{center}

Starting from the companion basis $\Pi$ for $\overrightarrow{D}_n$ and applying the corresponding sequence of outward companion basis mutations, it is routinely checked that we obtain the following companion basis for $\Xi_m$:

\begin{eqnarray*}
\label{companion basis for delta}
\beta_x^{\Xi_m} & = & \alpha_x \qquad \textrm{for all $x \neq n-1,n$};\\
\beta_{n-1}^{\Xi_m} & = & s_{\alpha_{m+1}}s_{\alpha_{m+2}}\cdots s_{\alpha_{n-2}} (\alpha_{n-1}) = \alpha_{m+1} + \alpha_{m+2} + \ldots + \alpha_{n-2} + \alpha_{n-1};\\
\beta_n^{\Xi_m} & = & s_{\alpha_{m+1}}s_{\alpha_{m+2}}\cdots s_{\alpha_{n-2}} (\alpha_{n}) = \alpha_{m+1} + \alpha_{m+2} + \ldots + \alpha_{n-2} + \alpha_{n}.
\end{eqnarray*}

Choose an ordered pair of distinct (if possible) end vertices in $\Gamma^{(1)}$ with $v_1$ being the second vertex in this pair. Applying our labelling procedure from the previous section, we then obtain a labelling of $\Gamma^{(1)}$. Apply the same labels to the corresponding vertices of $\Gamma$ (noting that the vertex $v_1$ gets the label $m$). Likewise, choose an ordered pair of distinct (if possible) end vertices in $\Gamma^{(2)}$ with $v_2$ being the first vertex in this pair. Applying our labelling procedure, we obtain a labelling of $\Gamma^{(2)}$ with the labels $\{ 1, \ldots , n-m-2 \}$. Add $m$ to each label, giving a labelling of $\Gamma^{(2)}$ with the labels $\{m+1, \ldots , n-2 \}$, and apply the same labels to the corresponding vertices of $\Gamma$ (noting that the vertex $v_2$ gets the label $m+1$). Finally, label the remaining two vertices of $\Gamma$ with $n-1$ and $n$, as shown in the diagram below:

 \begin{center}
 \psfragscanon
 \psfrag{n-1}{\mbox{\small $n-1$}}
 \psfrag{n}{\mbox{\small $n$}}
 \psfrag{m}{\mbox{\small $m$}}
 \psfrag{m+1}{\mbox{\small $m+1$}}
 \psfrag{g1}{\mbox{\small $\Gamma^{(1)}$}}
 \psfrag{g2}{\mbox{\small $\Gamma^{(2)}$}}
 \psfrag{g}{$\Gamma$:}
 \includegraphics[scale=.50]{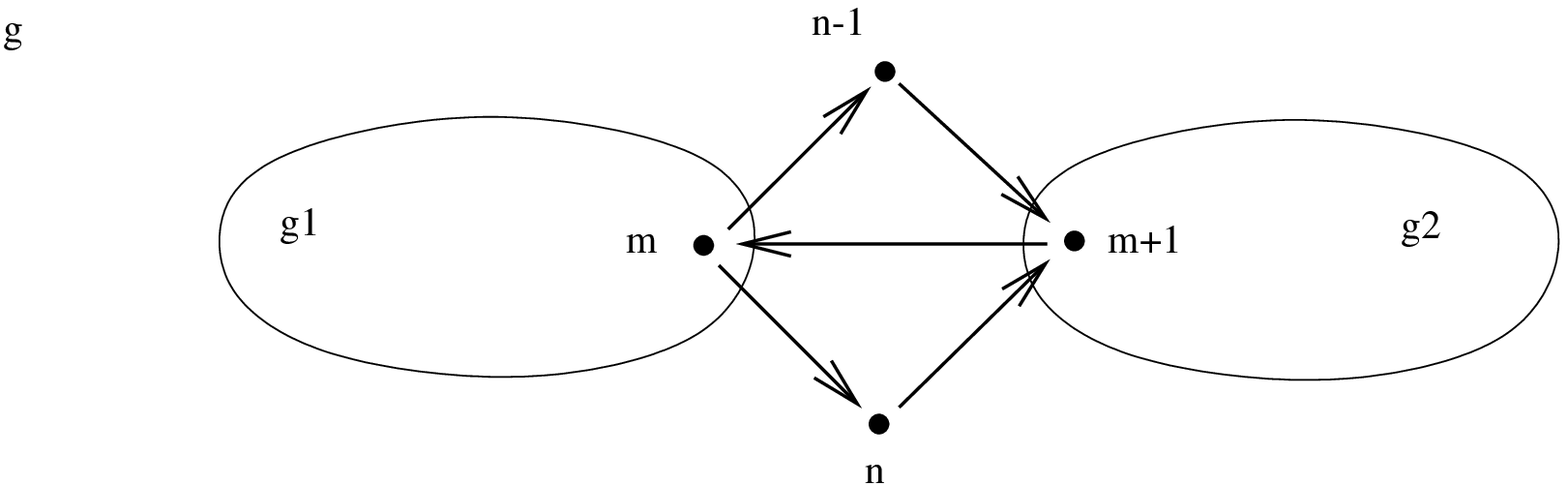}
 \end{center}

We saw in the previous section that by applying quiver mutations to the labelled quiver $\overrightarrow{A}_m$, without mutations at $1$ or $m$, we can obtain a labelled quiver $\widetilde{\Gamma}^{(1)}$ quasi-isomorphic to $\Gamma^{(1)}$. In addition, by applying quiver mutations to the labelled quiver

\vspace{6pt}

 \begin{center}
 \psfragscanon
 \psfrag{m+1}{\mbox{\small $m+1$}}
 \psfrag{m+2}{\mbox{\small $m+2$}}
 \psfrag{n-3}{\mbox{\small $n-3$}}
 \psfrag{n-2}{\mbox{\small $n-2$}}
 \psfrag{a}{$\overrightarrow{A}_{n-m-2}^{+m}$:}
 \includegraphics[scale=.50]{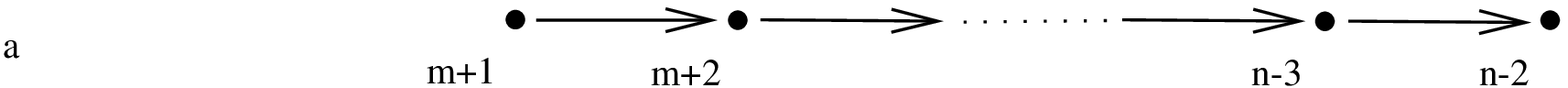}
 \end{center}

\noindent without mutations at $m+1$ or $n-2$, we can obtain a labelled quiver $\widetilde{\Gamma}^{(2)}$ quasi-isomorphic to $\Gamma^{(2)}$. Since mutation at $m$ or $m+1$ is not required, performing the corresponding sequences of mutations at the vertices of the quiver $\Xi_m$, we obtain a labelled quiver $\widetilde{\Gamma}$ quasi-isomorphic to $\Gamma$. (Note that here we are identifying the quiver $\overrightarrow{A}_m$ with the full subquiver of $\Xi_m$ on the vertices $\{1, \ldots , m\}$ and we are identifying the quiver $\overrightarrow{A}_{n-m-2}^{+m}$ with the full subquiver of $\Xi_m$ on the vertices $\{m+1, \ldots , n-2\}$.)

\bthm
\label{type D type II}
Let $\Gamma$ be a type II quiver of mutation type $D_n$, labelled as specified above. The set $\{ \beta_i \colon i \in \Gamma_0 \} \subseteq \Phi$ given by
$$
\beta_i = \left \{
\begin{array}{cl}
\alpha_i + \ldots + \alpha_{i+a} & \textrm{if } i \textrm{ is a primary vertex and } i+a \textrm{ is the}\\
 & \textrm{complementary vertex corresponding to } i,\\
\alpha_{m+1} + \alpha_{m+2} + \ldots + \alpha_{n-2} + \alpha_{n-1} & \textrm{if } i = n-1,\\
\alpha_{m+1} + \alpha_{m+2} + \ldots + \alpha_{n-2} + \alpha_{n} & \textrm{if } i = n,\\
\alpha_i & \textrm{otherwise}
\end{array}
\right.
$$
is a companion basis for $\Gamma$.
\ethm

\begin{proof}
Consider the sequence of quiver mutations taking us from $\Xi_m$ to $\widetilde{\Gamma}$. Starting from our given companion basis for $\Xi_m$, we may obtain a companion basis for $\widetilde{\Gamma}$ (and hence also for $\Gamma$) by performing the corresponding sequence of inward companion basis mutations.

Since we have $\beta_i^{\Xi_m} = \alpha_i$ for $i = 1, \ldots , m$ and $i = m+1, \ldots , n-2$, it follows similarly to the proof of Theorem~\ref{type D type I} that the companion basis elements we obtain associated to the vertices of $\widetilde{\Gamma}$ belonging to $\widetilde{\Gamma}^{(1)}$ and $\widetilde{\Gamma}^{(2)}$ match those given in the statement of the theorem. Moreover, since no companion basis mutations are performed at the vertices $m$ or $m+1$, it follows immediately from Theorem~\ref{companion basis mutation theorem} that the companion basis elements we obtain associated to the vertices $n-1$ and $n$ are the same as those appearing in our companion basis for $\Xi_m$.
\end{proof}

\underline{Type III}: We have that $\Gamma = \Gamma^{(0)} \underline{\sqcup} \left( \Gamma^{(1)}, \Gamma^{(2)} \right)$, with the end vertex $v_1$ of $\Gamma^{(1)}$ and the end vertex $v_2$ of $\Gamma^{(2)}$ identified respectively with the vertices $c_1$ and $c_2$ of the skeleton quiver $\Gamma^{(0)}$, as shown previously. Here, $\Gamma^{(1)}$ is a quiver of mutation type $A_m$ and $\Gamma^{(2)}$ is a quiver of mutation type $A_{n-m-2}$ (for some $1 \leq m \leq n-3$).

By proceeding as in the Type II case above, we obtain labellings of the vertices of $\Gamma^{(1)}$ and $\Gamma^{(2)}$ which together induce a labelling of the vertices of $\Gamma$, after labelling the remaining two vertices of $\Gamma$ with $n-1$ and $n$, as illustrated below:

 \begin{center}
 \psfragscanon
 \psfrag{n-1}{\mbox{\small $n-1$}}
 \psfrag{n}{\mbox{\small $n$}}
 \psfrag{m}{\mbox{\small $m$}}
 \psfrag{m+1}{\mbox{\small $m+1$}}
 \psfrag{g1}{\mbox{\small $\Gamma^{(1)}$}}
 \psfrag{g2}{\mbox{\small $\Gamma^{(2)}$}}
 \psfrag{g}{$\Gamma$:}
 \includegraphics[scale=.50]{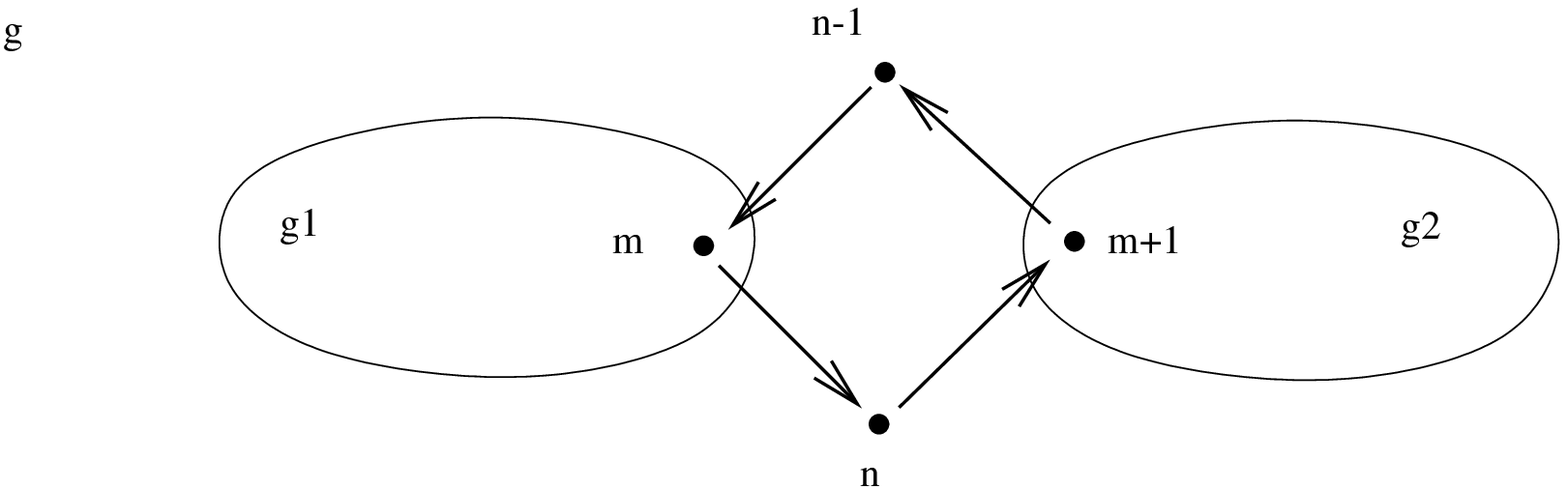}
 \end{center}

\bthm
\label{type D type III}
Let $\Gamma$ be a type III quiver of mutation type $D_n$, labelled as specified above. The set $\{ \beta_i \colon i \in \Gamma_0 \} \subseteq \Phi$ given by
$$
\beta_i = \left \{
\begin{array}{cl}
\alpha_i + \ldots + \alpha_{i+a} & \textrm{if } i \textrm{ is a primary vertex and } i+a \textrm{ is the}\\
 & \textrm{complementary vertex corresponding to } i,\\
\alpha_{m} + \ldots + \alpha_{n-2} + \alpha_{n-1} & \textrm{if } i = m,\\
\alpha_{m+1} + \ldots + \alpha_{n-2} + \alpha_{n-1} & \textrm{if } i = n-1,\\
\alpha_{m+1} + \ldots + \alpha_{n-2} + \alpha_{n} & \textrm{if } i = n,\\
\alpha_i & \textrm{otherwise}
\end{array}
\right.
$$
is a companion basis for $\Gamma$.
\ethm

\begin{proof}
Consider the quiver $\Gamma'$ obtained from $\Gamma$ by applying a quiver mutation at the vertex $n-1$:

 \begin{center}
 \psfragscanon
 \psfrag{n-1}{\mbox{\small $n-1$}}
 \psfrag{n}{\mbox{\small $n$}}
 \psfrag{m}{\mbox{\small $m$}}
 \psfrag{m+1}{\mbox{\small $m+1$}}
 \psfrag{g1}{\mbox{\small $\Gamma^{(1)}$}}
 \psfrag{g2}{\mbox{\small $\Gamma^{(2)}$}}
 \psfrag{g}{$\Gamma'$:}
 \includegraphics[scale=.50]{typeIIlabelled.eps}
 \end{center}

Theorem~\ref{type D type II} provides us with a companion basis for $\Gamma'$. Mutating the quiver $\Gamma'$ at the vertex $n-1$ takes us back to our original quiver $\Gamma$. Applying the corresponding inward companion basis mutation at the vertex $n-1$ to our companion basis for $\Gamma'$ then establishes the result. We note that only the companion basis element associated to the vertex $m$ changes, as the only arrow going into $n-1$ has source $m$, and furthermore,
$$
\beta_m^{\Gamma} = s_{\beta_{n-1}^{\Gamma'}}(\beta_m^{\Gamma'}) = s_{\alpha_{m+1} + \ldots + \alpha_{n-2} + \alpha_{n-1}}(\alpha_{m}) = \alpha_{m} + \ldots + \alpha_{n-2} + \alpha_{n-1}.
$$
\end{proof}

\underline{Type IV}: Either we have that $\Gamma$ is an oriented $n$-cycle or we have that $\Gamma = \Gamma^{(0)} \underline{\sqcup} \left( \Gamma^{(1)}, \ldots, \Gamma^{(r)} \right)$, with the end vertex $v_i$ of $\Gamma^{(i)}$ identified with the vertex $c_i$ of $\Gamma^{(0)}$ for each $i$, $1 \leq i \leq r$. Consider the latter case further, supposing that $\Gamma^{(i)}$ has $a_i$ vertices for each $i \neq 0$, so that $\Gamma^{(i)}$ is a quiver of mutation type $A_{a_i}$. We recall that $\Gamma^{(0)}$ consists of a central cycle decorated with $r$ spikes. The central cycle is an oriented $m$ cycle, for some $m\geq 3$. The $i$th spike, which we denote $S_i$, is an oriented 3-cycle $p_i \stackrel{\alpha_i}{\longrightarrow} s_i \longrightarrow c_i \longrightarrow p_i$ containing both the arrow $p_i \stackrel{\alpha_i}{\longrightarrow} s_i$ and the vertex $c_i$. The arrows $\alpha_1, \ldots , \alpha_r$ form a subset of the arrows of the central cycle, and appear in clockwise order around this cycle. The vertices $p_1, \ldots , p_r$ (and likewise $s_1,\ldots,s_r$) are therefore distinct and appear in clockwise order around the central cycle. It can happen that $s_i = p_{i+1}$ for some values of $i$, and we may also have $s_r = p_1$. It will be useful to note that, due to the vertex identifications, $n = m+ \sum_{i = 1}^{r}{a_i}$.

We first apply a labelling to $\Gamma$. It will be convenient for us to consider two cases separately. The first case being that where at least one arrow in the central cycle doesn't belong to a spike (that is, where $r<m$), and the second where every arrow on the central cycle belongs to a spike (that is, where $r=m$). We note that, by definition, we must have $r \leq m$.

Case 1 ($0 \leq r<m$): We initially suppose that $0<r<m$. By relabelling if necessary, we can choose $S_1$ such that $p_1 \neq s_r$. (That is, we can choose $S_1$ such that the 1st and the $r$th spikes do not meet in a vertex.)

We start by labelling the vertices of the central cycle. It will be helpful to use $d_i$ to denote the number of (non-spiked) arrows between $S_i$ and $S_{i+1}$ on the central cycle, for $1 \leq i \leq r-1$. In addition, we will use $d_r$ to denote the number of arrows between $S_r$ and $S_1$ on the central cycle. We label the vertex $p_1$ with 1 and then, travelling in the clockwise direction around the central cycle, label successive vertices in increments of 1, except that for each vertex $p_i$, $1\leq i \leq r$, the subsequent vertex $s_i$ is labelled $p_i+a_i+1$ (note that $p_i$ in this expression refers to the label of the vertex $p_i$). So, in general, for $1\leq i \leq r$, the vertex $p_i$ will be labelled $\sum_{j=1}^{i-1} a_j + \sum_{j=1}^{i-1} d_j + i$ and the vertex $s_i$ will be labelled $p_i + a_i + 1 = \sum_{j=1}^{i} a_j + \sum_{j=1}^{i-1} d_j + i + 1$. (It can happen that $s_i = p_{i+1}$ for some values of $i\neq r$. The expressions for $s_i$ and $p_{i+1}$ agree in this case, since $d_i = 0$.) By abuse of notation, we will identify vertex names with their labels, and use these interchangeably.

We complete the labelling of $\Gamma$ by, for each $i$, $1\leq i \leq r$, assigning the $a_i$ labels (strictly) in-between $p_i$ and $s_i$ to the vertices of the quiver $\Gamma^{(i)}$. We choose an ordered pair of distinct (if possible) end vertices in $\Gamma^{(i)}$, with $v_i$ being the second vertex in this pair. Applying our labelling procedure we obtain a labelling of $\Gamma^{(i)}$ with the labels $\{ 1, \ldots, a_i \}$. We then add $p_i$ to each label, giving a labelling of $\Gamma^{(i)}$ with the labels $\{ p_i+1, \ldots, p_i+a_i \}$, and assign these labels to the corresponding vertices of $\Gamma$.

If $r=0$ (that is, if $\Gamma^{(0)}$ has no spikes, so that $\Gamma = \Gamma^{(0)}$ is an oriented $n$-cycle), we label $\Gamma$ as follows: We first choose an arbitrary vertex of $\Gamma$, which we label with 1. Then, travelling in the clockwise direction around $\Gamma$, we label successive vertices in increments of 1 until the labelling is completed.

Case 2 ($r=m$): In this case, we have that every arrow on the central cycle belongs to a spike. So, $s_r$ and $p_1$ are the same vertex, and also, $s_i = p_{i+1}$ for all $i$, $1\leq i \leq r-1$. In addition, $d_i = 0$ for all $i$, $1\leq i \leq r$. Again, we start by labelling the vertices of the central cycle. These vertices are precisely $p_1, \ldots, p_r$. For each $i$, $1 \leq i \leq r$, we label $p_i$ with $\sum_{j=1}^{i-1} a_j + i$. In particular, $p_1$ is the first vertex to be labelled, and gets the label 1.

We label the vertices of $\Gamma^{(1)}, \ldots, \Gamma^{(r-1)}$ precisely as in Case 1 above (noting that $s_i = p_{i+1}$ here, for $1\leq i \leq r-1$). It just remains to label the vertices of $\Gamma^{(r)}$. We choose an ordered pair (distinct if possible) of end vertices in $\Gamma^{(r)}$, with $v_r$ being the second vertex in this pair. Applying our labelling procedure we obtain a labelling of $\Gamma^{(r)}$ with the labels $\{ 1, \ldots, a_r \}$. We then add $p_r$ to each label, giving a labelling of $\Gamma^{(r)}$ with the labels $\{ p_r+1, \ldots, p_r+a_r \}$. Note here that $v_r$ is the final vertex to be labelled and gets the label $p_r + a_r$, and moreover, that $p_r + a_r = n$ (this can be seen by referring to expressions for $n$ and $p_r$ given previously, and using that $r = m$). Finally, we assign these labels to the corresponding vertices of $\Gamma$, thus completing the labelling of $\Gamma$.

We immediately state the results providing a companion basis for $\Gamma$. The first of these applies in Case 1 and the second in Case 2.

\bthm
\label{type D type IV case 1}
Let $\Gamma$ be a type IV quiver of mutation type $D_n$ ($0 \leq r<m$ case), labelled as specified above. The set $\{ \beta_i \colon i \in \Gamma_0 \} \subseteq \Phi$ given by
$$
\beta_i = \left \{
\begin{array}{cl}
\alpha_i + \ldots + \alpha_{i+a} & \textrm{if } i \textrm{ is a primary vertex and } i+a \textrm{ is the}\\
 & \textrm{complementary vertex corresponding to } i,\\
\alpha_{p_j} + \ldots + \alpha_{p_j+a_j} & \textrm{if } i = p_j \textrm{ for some } 1\leq j \leq r,\\
\alpha_{1} + \ldots + \alpha_{n-2} + \alpha_{n} & \textrm{if } i = n,\\
\alpha_i & \textrm{otherwise}
\end{array}
\right.
$$
is a companion basis for $\Gamma$. (Noting that $p_j + a_j$ is necessarily the label of the vertex $v_j$ helps with reading-off this companion basis from the labelled quiver $\Gamma$.)
\ethm

\bthm
\label{type D type IV case 2}
Let $\Gamma$ be a type IV quiver of mutation type $D_n$ ($r=m$ case), labelled as specified above. The set $\{ \beta_i \colon i \in \Gamma_0 \} \subseteq \Phi$ given by
$$
\beta_i = \left \{
\begin{array}{cl}
\alpha_i + \ldots + \alpha_{i+a} & \textrm{if } i \textrm{ is a primary vertex and } i+a \textrm{ is the}\\
 & \textrm{complementary vertex corresponding to } i,\\
\alpha_{p_j} + \ldots + \alpha_{p_j+a_j} & \textrm{if } i = p_j \textrm{ for some } 1\leq j \leq r-1,\\
\alpha_1 + \ldots +\alpha_{p_r-1} + 2\alpha_{p_r} + \ldots + 2\alpha_{n-2} + \alpha_{n-1} + \alpha_n & \textrm{if } i = p_r,\\
\alpha_{1} + \ldots + \alpha_{n-2} + \alpha_{n} & \textrm{if } i = n,\\
\alpha_i & \textrm{otherwise}
\end{array}
\right.
$$
is a companion basis for $\Gamma$. (Noting that $p_j + a_j$ is necessarily the label of the vertex $v_j$ helps with reading-off this companion basis from the labelled quiver $\Gamma$.)
\ethm

We devote the remainder of this section to completing the proofs of Theorems~\ref{type D type IV case 1} and \ref{type D type IV case 2}. The first step, in each case, is to show how we can mutate the quiver $\overrightarrow{D}_n$ into a quiver $\widetilde{\Gamma}$ quasi-isomorphic to $\Gamma$. This will enable us to construct a companion basis for $\Gamma$ by simultaneously performing the corresponding companion basis mutations, starting from the companion basis $\Pi$ for $\overrightarrow{D}_n$. Much of the procedure for transforming $\overrightarrow{D}_n$ into a quiver quasi-isomorphic to $\Gamma$ is common to both cases, so we only consider the two cases separately at the point at which it becomes necessary.

The first step is to mutate $\overrightarrow{D}_n$ into an oriented $n$-cycle. Starting from $\overrightarrow{D}_n$, applying consecutive quiver mutations at the vertices $n-1, n-2, \ldots, 1$, we obtain the oriented labelled $n$-cycle $N_n$, as shown below:

 \begin{center}
 \psfragscanon
 \psfrag{1}{\mbox{\small $1$}}
 \psfrag{2}{\mbox{\small $2$}}
 \psfrag{3}{\mbox{\small $3$}}
 \psfrag{4}{\mbox{\small $4$}}
 \psfrag{5}{\mbox{\small $5$}}
 \psfrag{n-1}{\mbox{\small $n-1$}}
 \psfrag{n}{\mbox{\small $n$}}
 \psfrag{N}{$N_n$:}
 \includegraphics[scale=.50]{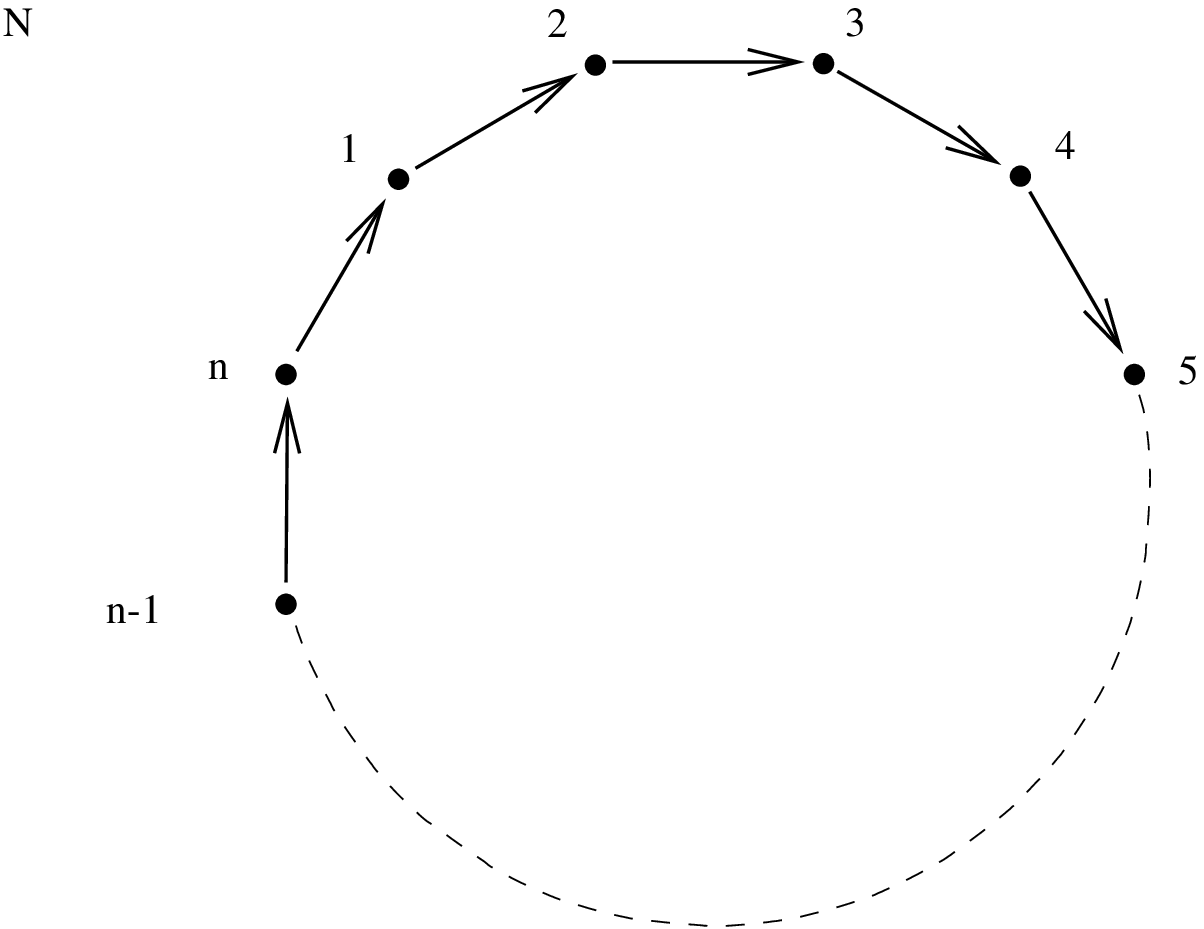}
 \end{center}

Starting from the companion basis $\Pi$ for $\overrightarrow{D}_n$ and applying the corresponding sequence of outward companion basis mutations, we obtain the following companion basis for $N_n$:

\begin{eqnarray*}
\beta_x^{N_n} & = & \alpha_x \qquad \textrm{for all } x \neq n; \\
\beta_n^{N_n} & = & s_{\alpha_1} \cdots s_{\alpha_{n-2}} (\alpha_n) = \alpha_1 + \cdots + \alpha_{n-2} + \alpha_n.
\end{eqnarray*}

We note that this proves Theorem~\ref{type D type IV case 1} in the special case where $r=0$. It remains to deal with the cases where $r>0$. That is, the cases where $\Gamma$ has spikes. The next step is to mutate $N_n$ to produce a labelled quiver $\Theta$ which has (the labelled quiver) $\Gamma^{(0)}$ as a full subquiver, and such that $\Theta$ is obtained from the disjoint union $\Gamma^{(0)} \sqcup A_{a_1}^{+p_1} \sqcup \cdots \sqcup A_{a_r}^{+p_r}$ by identifying the vertex $p_i+a_i$ with the vertex $c_i$ for each $i$, $1 \leq i \leq r$. From this intermediate step, we will then be able construct $\widetilde{\Gamma}$ by repeatedly applying Proposition~\ref{labelling by mutation}.

\begin{figure}[ht]
 \begin{center}
 \psfragscanon
 \psfrag{p1}{\mbox{\tiny $p_1 = 1$}}
 \psfrag{p1+1}{\mbox{\tiny $p_1+1$}}
 \psfrag{p1+a1}{\mbox{\tiny $p_1+a_1$}}
 \psfrag{p2}{\mbox{\tiny $p_2$}}
 \psfrag{p2+1}{\mbox{\tiny $p_2+1$}}
 \psfrag{p2+a2}{\mbox{\tiny $p_2+a_2$}}
 \psfrag{p3}{\mbox{\tiny $p_3$}}
 \psfrag{p3+1}{\mbox{\tiny $p_3+1$}}
 \psfrag{p3+a3}{\mbox{\tiny $p_3+a_3$}}
 \psfrag{pi}{\mbox{\tiny $p_i$}}
 \psfrag{pi+1}{\mbox{\tiny $p_i+1$}}
 \psfrag{pi+ai}{\mbox{\tiny $p_i+a_i$}}
 \psfrag{pr}{\mbox{\tiny $p_r$}}
 \psfrag{pr+1}{\mbox{\tiny $p_r+1$}}
 \psfrag{pr+ar}{\mbox{\tiny $p_r+a_r$}}
 \psfrag{T}{$\Theta$:}
 \includegraphics[scale=.50]{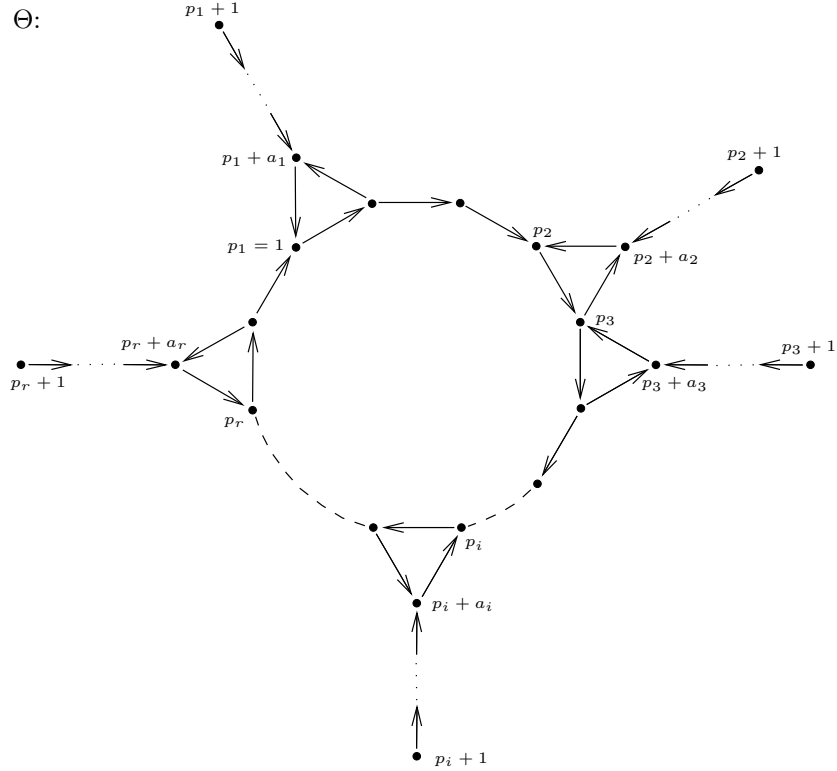}
 \caption{Diagram of $\Theta$ in Case 1.}
 \label{theta in case 1}
 \end{center}
\end{figure}

\begin{figure}[ht]
 \begin{center}
 \psfragscanon
 \psfrag{p1}{\mbox{\tiny $p_1 = 1$}}
 \psfrag{p1+1}{\mbox{\tiny $p_1+1$}}
 \psfrag{p1+a1}{\mbox{\tiny $p_1+a_1$}}
 \psfrag{p2}{\mbox{\tiny $p_2$}}
 \psfrag{p2+1}{\mbox{\tiny $p_2+1$}}
 \psfrag{p2+a2}{\mbox{\tiny $p_2+a_2$}}
 \psfrag{p3}{\mbox{\tiny $p_3$}}
 \psfrag{p3+1}{\mbox{\tiny $p_3+1$}}
 \psfrag{p3+a3}{\mbox{\tiny $p_3+a_3$}}
 \psfrag{pi}{\mbox{\tiny $p_i$}}
 \psfrag{pi+1}{\mbox{\tiny $p_i+1$}}
 \psfrag{pi+ai}{\mbox{\tiny $p_i+a_i$}}
 \psfrag{pr}{\mbox{\tiny $p_r$}}
 \psfrag{pr+1}{\mbox{\tiny $p_r+1$}}
 \psfrag{pr+ar}{\mbox{\tiny $p_r+a_r$}}
 \psfrag{T}{$\Theta$:}
 \includegraphics[scale=.50]{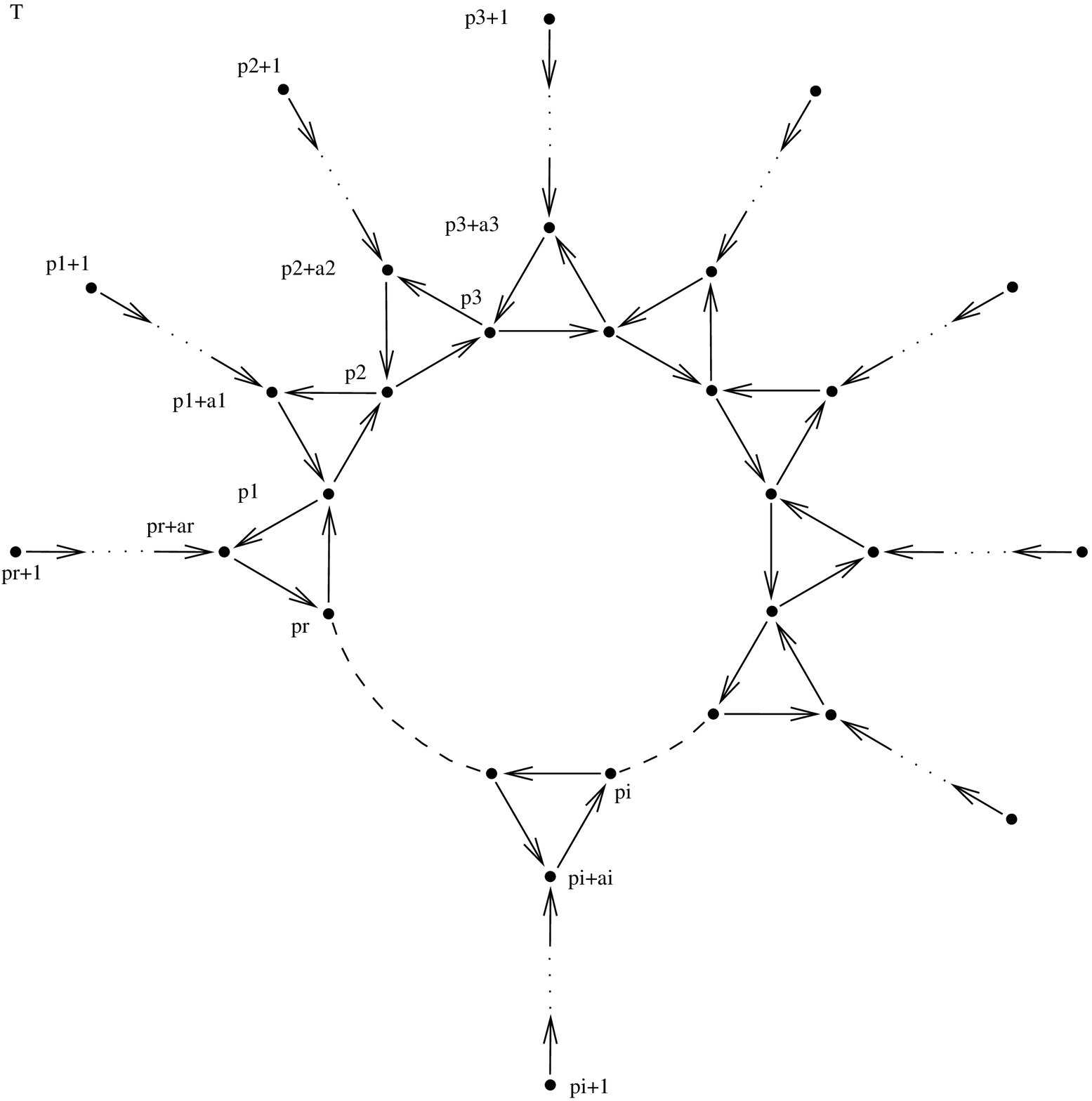}
 \caption{Diagram of $\Theta$ in Case 2.}
 \label{theta in case 2}
 \end{center}
\end{figure}

We produce $\Theta$ by starting from $N_n$ and mutating successively at the vertices $p_1 +1, \ldots , p_1 + a_1$ followed by $p_2 +1, \ldots, p_2 + a_2, \ldots, p_i +1, \ldots, p_i + a_i, \ldots$, and then finally $p_r +1, \ldots, p_r + a_r$. Figures~\ref{theta in case 1} and \ref{theta in case 2} illustrate $\Theta$ in cases 1 and 2 respectively.

Case 1: Performing the corresponding sequence of inward companion basis mutations to the above companion basis for $N_n$, it is easily observed (as in the proof of Proposition~\ref{companion basis for quasi-isomorphic quiver}) that we obtain the following companion basis for $\Theta$:

\begin{eqnarray*}
\beta_{p_i}^{\Theta} & = & s_{\beta_{p_i + a_i}^{N_n}} \cdots s_{\beta_{p_i + 1}^{N_n}} (\beta_{p_i}^{N_n}) \\
 & = & s_{\alpha_{p_i + a_i}} \cdots s_{\alpha_{p_i +1}} (\alpha_{p_i}) \\
 & = & \alpha_{p_i} + \ldots + \alpha_{p_i + a_i} \qquad \textrm{for } 1 \leq i \leq r; \\
\beta_n^{\Theta} & = & \beta_n^{N_n} \\
 & = & \alpha_1 + \cdots + \alpha_{n-2} + \alpha_n; \\
\beta_i^{\Theta} & = & \beta_i^{N_n} \\
 & = & \alpha_i \qquad \textrm{otherwise}.
\end{eqnarray*}

Case 2: The situation here is identical to that for Case 1, except that now we have $p_r + a_r = n$ and so $\beta_{p_r}^{\Theta}$ is instead given by:

\begin{eqnarray*}
\beta_{p_r}^{\Theta} & = & s_{\beta_{p_r + a_r}^{N_n}} \cdots s_{\beta_{p_r + 1}^{N_n}} (\beta_{p_r}^{N_n}) \\
 & = & s_{\beta_{n}^{N_n}} s_{\beta_{n-1}^{N_n}} \cdots s_{\beta_{p_r + 1}^{N_n}} (\beta_{p_r}^{N_n}) \\
 & = & s_{\alpha_1 + \cdots + \alpha_{n-2} + \alpha_n} s_{\alpha_{n-1}} \cdots s_{\alpha_{p_r +1}} (\alpha_{p_r}) \\
 & = & s_{\alpha_1 + \cdots + \alpha_{n-2} + \alpha_n} (\alpha_{p_r} + \cdots + \alpha_{n-1}) \\
 & = & \alpha_1 + \ldots +\alpha_{p_r-1} + 2\alpha_{p_r} + \ldots + 2\alpha_{n-2} + \alpha_{n-1} + \alpha_n.
\end{eqnarray*}

The procedure for obtaining a quiver quasi-isomorphic to $\Gamma$ from $\Theta$ is the same in both cases: For each $i$, we can mutate $A_{a_i}^{+p_i}$ into a labelled quiver $\widetilde{\Gamma}^{(i)}$ quasi-isomorphic to $\Gamma^{(i)}$, without mutation at the vertex $p_i+a_i$, as in Proposition~\ref{labelling by mutation}. The $i$th step of the procedure, $1\leq i \leq r$, is to perform the mutations taking us from $A_{a_i}^{+p_i}$ to $\widetilde{\Gamma}^{(i)}$ to the corresponding vertices of the labelled quiver $\Theta^{(i-1)}$, thereby obtaining a new labelled quiver which we denote $\Theta^{(i)}$ (we note that as a starting point, $\Theta^{(0)}$ is used to denote the quiver $\Theta$). Upon completion of the $r$th step of the procedure, it is easily confirmed that we obtain a labelled quiver $\widetilde{\Gamma} (= \Theta^{(r)})$ quasi-isomorphic to $\Gamma$.

We are now ready to complete the proofs of Theorems~\ref{type D type IV case 1} and \ref{type D type IV case 2}. In each case, we essentially obtain a companion basis for $\widetilde{\Gamma}$ (and hence also for $\Gamma$) by applying the sequence of inward companion basis mutations corresponding to the quiver mutations performed in transforming $\Theta$ into $\widetilde{\Gamma}$. However, much of the detail is hidden due to the fact that we are able to take advantage of our companion basis construction procedure for Dynkin type $A$.

\begin{proof}[Proof of Theorem~\ref{type D type IV case 1}]
\label{type D type IV case 1 proof}
For each $i$, $1\leq i \leq r$, we have that $\beta_j^{\Theta} = \alpha_j$ for all vertices $j$ of $A_{a_i}^{+p_i}$. It therefore follows similarly to the proof of Theorem~\ref{type D type I} that the result holds for the vertices of $\widetilde{\Gamma}$ belonging to any $\widetilde{\Gamma}^{(i)}$, $i \neq 0$. That the result holds for the vertices of the central cycle of $\Gamma^{(0)}$ follows from the fact that no mutations are performed at the vertices $c_1, \ldots , c_r$ (those labelled $p_1+a_1, \ldots , p_r+a_r$) in the sequence of mutations taking us from $\Theta$ to $\widetilde{\Gamma}$. In particular, we must have $\beta_i^{\widetilde{\Gamma}} = \beta_i^{\Theta}$ for all vertices $i$ of the central cycle of $\Gamma^{(0)}$ (upon performing the corresponding sequence of inward companion basis mutations).
\end{proof}

\begin{proof}[Proof of Theorem~\ref{type D type IV case 2}]
\label{type D type IV case 2 proof}
That the result holds on the vertices of $\widetilde{\Gamma}$ belonging to $\widetilde{\Gamma}^{(1)}, \ldots , \widetilde{\Gamma}^{(r-1)}$ follows as in the proof of Theorem~\ref{type D type IV case 1} above. We now consider the full subquiver $(A_{a_r}^{+p_r})$ of $\Theta^{(r-1)}$ on the vertices $p_r+1, \ldots , p_r+a_r = n$:

 \begin{center}
 \psfragscanon
 \psfrag{pr+1}{\mbox{\small $p_r+1$}}
 \psfrag{pr+2}{\mbox{\small $p_r+2$}}
 \psfrag{pr+ar}{\mbox{\small $p_r+a_r = n$}}
 \psfrag{a}{$A_{a_r}^{+p_r}$:}
 \includegraphics[scale=.50]{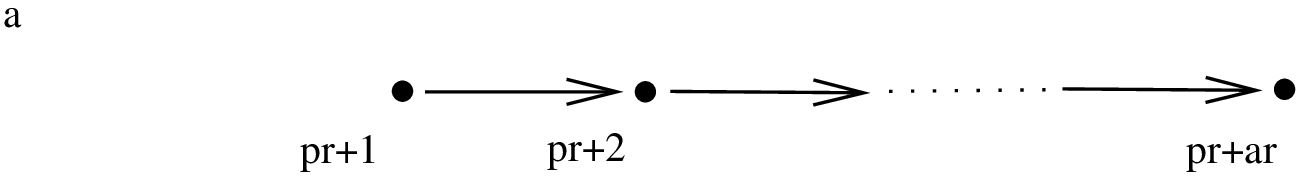}
 \end{center}

\noindent In the companion basis for $\Theta^{(r-1)}$ obtained from that for $\Theta$ by applying the inward companion basis mutations corresponding to the sequence of quiver mutations used to get from $\Theta$ to $\Theta^{(r-1)}$, we have that $\beta_i^{\Theta^{(r-1)}} = \alpha_i$ for $i = p_r+1, \ldots , n-1$ and $\beta_n^{\Theta^{(r-1)}} = \alpha_{1} + \ldots + \alpha_{n-2} + \alpha_{n}$. The proof of Proposition~\ref{labelling by mutation} provides us with the sequence of quiver mutations used to transform the labelled quiver $A_{a_r}^{+p_r}$ into the labelled quiver $\widetilde{\Gamma}^{(r)}$. In this sequence of mutations, it is clear that there is never an arrow going from the vertex $n$ into the vertex of mutation. Performing the corresponding sequence of inward companion basis mutations inside $\Theta^{(r-1)}$, it then follows in a similar manner to the proof of Proposition~\ref{companion basis for quasi-isomorphic quiver} that the result also holds for the vertices of $\widetilde{\Gamma}$ belonging to $\widetilde{\Gamma}^{(r)}$. Finally, we have that the result holds for the vertices of the central cycle of $\Gamma^{(0)}$, with similar reasoning to that given in the proof of Theorem~\ref{type D type IV case 1} above.
\end{proof}



\noindent \textbf{Acknowledgements:}
This work was supported by the Austrian Science Fund (FWF): Project Number P25141-N26. The author wishes to thank Robert Marsh for his continuing invaluable support. He also wishes to thank Karin Baur for her useful comments during the preparation of this article. The companion basis construction procedure for quivers of mutation type $A$ presented in this article previously appeared in the author's Ph.D. thesis \cite{Par1} (the justification given for it, however, is new). The author's Ph.D. was completed at the University of Leicester, under the supervision of Robert Marsh, and funded by the Engineering and Physical Sciences Research Council.

\end{document}